\documentclass[11pt,twoside]{article}
\usepackage{amsmath, amsthm, amsfonts, amssymb}
\usepackage{mathrsfs}
\usepackage[misc]{ifsym}
\usepackage[colorlinks,linkcolor=blue,anchorcolor=blue,citecolor=blue,urlcolor=black]{hyperref}
\usepackage{graphicx}
\usepackage{subfigure}
\usepackage{tikz}
\usepackage{enumerate}
\usepackage{anysize}
\usepackage{setspace}
\usepackage{float}
\usepackage{color}
\usepackage{chngcntr}
\everymath{\displaystyle}
\allowdisplaybreaks

\newfam\msbfam
\font\tenmsb=msbm5    \textfont\msbfam=\tenmsb \font\sevenmsb=msbm5
\scriptfont\msbfam=\sevenmsb \font\fivemsb=msbm5
\scriptscriptfont\msbfam=\fivemsb

\newfam\bigfam
\font\tenbig=msbm5 scaled \magstep2   \textfont\bigfam=\tenbig
\font\sevenbig=msbm7 scaled \magstep2 \scriptfont\bigfam=\sevenbig
\font\fivebig=msbm5 scaled \magstep2
\scriptscriptfont\bigfam=\fivebig

\numberwithin{equation}{section}
\newtheorem{theorem}{Theorem}[section]
\newtheorem{lemma}{Lemma}[section]

\newtheorem{remark}{Remark}[section]

\newtheorem{definition}{Definition}[section]

\begin{document}
\title{\bf Multilinear Square Operators Meet New Weight Functions}
\author{\bf Chunliang Li,   Shuhui  Yang  and Yan Lin$^*$  }
\renewcommand{\thefootnote}{}
\date{}
\maketitle
\footnotetext{This work was partially supported by the National 
Natural Science Foundation of China No. 12071052.}
\footnotetext{$^*$Corresponding author, E-mail:linyan@cumtb.edu.cn}
\begin{minipage}{13.5cm}
{\bf Abstract}
\quad
Via the new weight $A_{\vec p}^{\theta }(\varphi )$, the 
authors introduce a new class of multilinear square operators. The boundedness 
on the weighted Lebesgue space and the weighted Morrey space is obtained, 
respectively. Our results include the known results of the standard 
multilinear square operator and the weight $ A_{\vec p} $. Moreover, the 
results in this article seem to be new even for one-linear case.
\medskip
\end{minipage}

\section{Introduction}\label{sec1}

\quad\quad The multilinear theory has been attracted more attentions since the 
pioneering work of Coifman and Meyer \cite{b6}. In recent years, as a rapidly 
developing field of harmonic analysis, the multilinear theory has attracted the 
research interest of many scholars, including Kenig and Stein \cite{b17}, 
Lerner and Ombrosi \cite{b19}, Christ and Journé \cite{b5}, Grafakos and 
Torres \cite{b10,b11} et al. In 1978,  Coifman and Meyer \cite{b7} introduced a 
class of multilinear operators (bilinear, one dimensional) as a 
multilinearization of Littlewood--Paley $ g $-function,
\begin{align*}
\mathscr{A}(a, f):=\int_{0}^{\infty}\left(f * \phi_{t}\right)\left(a * 
\Phi_{t}\right) \frac{m(t)}{t} dt,
\end {align*}
where $  m\left( t \right) \in {L^\infty }(\mathbb{R})$, $ a \in 
BMO(\mathbb{R}) $, $ \phi _t(x):=\phi (x /t)t^{-n} $, $ \Phi _t(x):=\Phi (x 
/t)t^{-n} $, $ \widehat \phi $ and $ \widehat \Phi $  have compact 
support with $ 0\notin $ supp $\widehat \Phi$. They studied the $ L^2 $ 
estimate of this operator by using the notion of Carleson measure. Since then, 
many authors researched multilinearized Littlewood--Paley operators. In 1982, 
Yabuta \cite{b37} studied $H^1$, $ BMO $ and $L^p$ estimates of $ \mathscr{A} $ 
by weakening the assumptions. In 2002, Sato and Yabuta \cite{b28} showed that $ 
T_{g} $ is bounded from $ L^{p_{1}  } \times \cdots \times L^{p_{m}  } $ to $ 
L^{p} $ for multilinearized Littlewood--Paley operators,
\begin{align*}
T_{g}(\vec{f})(x):=\int_{0}^{\infty} 
\prod_{i=1}^{m}\left(\left(\phi_{i}\right)_{t} * f_i\right)(x) \frac{d t}{t}.
\end {align*}

In 2015, Xue and Peng et al \cite{b35} studied the strong and weak
estimates of the multilinear Littlewood--Paley $g$-function with convolution 
type,
\begin{align*}
g(\vec{f})(x):=\left(\int_{0}^{\infty}\left|\frac{1}{t^{mn}}\int_{({\mathbb{R}^n}
)^m}\psi(\frac{y_1}{t},\ldots , \frac{y_m}{t}) \prod _{j=1}^{m} 
f_j(x-y_j)d{\vec y}\right|^{2} \frac{d t}{t}\right)^{\frac{1}{2} }
\end {align*}
and its commutator on weighted Lebesgue spaces. 
In 2015, He and Xue et al \cite{b14} established the boundedness of multilinear 
square operator $T$ (see \eqref{1.1}) with non-convolution type on Campanato 
spaces.  Xue and Yan \cite{b36} established that the multilinear square 
function $ T $ is bounded from $ L^{p_{1} }\left ( \omega _{1} \right )  \times 
\cdots  \times 
L^{p_{m} } \left ( \omega _{m} \right ) $ to $ L^{p} \left ( \nu _{\vec{\omega 
} }  \right )$ with each $ p_{i}> 1 $, and $ T $ 
is bounded from $ L^{p_{1} }\left ( \omega _{1}  \right )  \times \cdots  
\times L^{p_{m} } \left ( \omega _{m}  \right ) $ to $ L^{p,\infty } \left ( 
\nu _{\vec{\omega } } \right )$ when there is a $  p_{i}= 1 $. 
In 2015, Li and Song \cite{b20} proved the vector-valued weighted norm 
boundedness for the multilinear square operator $ T $ on weighted Lebesgue 
spaces. In 2016, Si and Xue \cite{b30} attenuated the classical kernel 
conditions in \cite{b36}. They studied the strong $ (L^{p_1}(\omega _1)\times 
\cdots \times L^{p_m}(\omega _m),L^{p}(\nu  _{\vec \omega }) ) $ and the weak $ 
(L^{p_1}(\omega _1)\times \cdots \times L^{p_m}(\omega _m),L^{p,\infty }(\nu  
_{\vec \omega }) ) $ type estimates for multilinear square functions with 
kernel of Dini’s type. In 2017, Si \cite{b29} established the weighted strong 
type and weighted end-point weak type estimates for the  iterated commutators $ 
T_{\prod \vec b,q} $, where the  iterated commutators $ T_{\prod \vec b,q} $ 
were generated by the multilinear square fucntions $ T $ and $ BMO $ functions. 
In 2018, Si and Xue \cite{b31} obtained that the iterated commutators $ 
T_{\prod \vec b} $ were bounded from product Lebesgue spaces into Lebesgue 
spaces, Lipschitz spaces, and Triebel–Lizorkin spaces, where $ T_{\prod \vec b} 
$ were generated by the multilinear square 
fucntions with Dini's type kernels and Lipschitz functions. In 2018, Hormozi, 
Si and Xue \cite{b15} obtained the boundedness of multilinear square functions  
with non-smooth kernels on weighted Lebesgue spaces. The study on the 
multilinear square functions have crucial applications in PDEs and other 
fields. For further details on the theory of multilinear square functions and 
their applications, one can see \cite{b3,b4,b9,b13} and related references.

Let a locally integrable function $K(y_0,y_1,\ldots,y_m)$ be defined away 
from the diagonal $y_0=y_1=\ldots=y_m$ in $(\mathbb{R}^n)^{m+1}$ and $ 
K_t(x,y_1,\ldots ,y_m):=(1/t^{mn})K(x/t,y_1/t,\ldots,y_m/t)$, for any $ t\in 
(0,\infty )$. The multilinear square function is given by
\begin{equation}\label{1.1}	
T(\vec{f})(x):={\left( {\int_0^\infty  {{{\left|{\int_{{{\left({\mathbb{R}^n} 
\right)}^m}} {{K_t}\left( {x,{y_1},\ldots,{y_m}} \right) \cdot 
\prod_{j=1}^{m}f_{j}\left(y_{j}\right)d{\vec y}} }
\right|}^2}\frac{{dt}}{t}} } \right)^{\frac{1}{2}}},
\end{equation}
whenever $x\notin  {\textstyle \bigcap_{j=1}^{m}}  {\rm supp } f_j$ 
and each $ f_{j} \in C_{c}^{\infty}\left(\mathbb{R}^{n}\right) $.

It is well known that weighted inequalities play a significant role in the 
harmonic analysis. In recent years, weighted inequalities has attracted great 
research interest from many scholars. In \cite{b8,b23,b27}, the authors studied 
the classical Muckenhoupt's class weight functions.

Let $\omega(x) \in L_{loc}^1(\mathbb{R}^n)$ and be nonnegative. We say that 
$\omega \in A_p$ for $1<p<\infty$, $ 1/p+1/p'=1 $ and all cubes, if and only if 
there exists a positive constant $C$ such that
\begin{align*}
\left(\frac{1}{|Q|}\int_Q\omega(y)dy\right)^{\frac{1}{p}}
\left(\frac{1}{|Q|}\int_Q\omega(y)^{\frac{-1}{p-1}}dy\right)^
{\frac{1}{p'}}\leq C.
\end {align*}
When $ p=1 $, $\omega \in A_1$ if and only if there exists a positive constant 
$C$ such that
\begin{align*}
\left(\frac{1}{|Q|}\int_Q\omega(y)dy\right)\leq C\mathop{\inf}\limits_{x\in 
{\rm Q}}\omega(x).
\end {align*}
Via the classical Muckenhoupt's class weight functions, 
numrous authors have studied the extrapolation of $ A_p$. In 2012, Tang 
\cite{b32} obtained new weighted norm 
inequalities for pseudo-differential operators with smooth symbols, where the 
weights are $A_p^{\infty}(\varphi)$. For some other studies of 
$A_p^{\infty}(\varphi)$, one can refer to \cite{b12,b16,b21,b33} and so on. 

We say that $\omega \in A_p^{\theta}(\varphi)$ for $ 1< p< \infty $, $ 
1/p+1/p'=1$ and $\theta\geq0$, if and only if there exists a 
positive constant $C$ such that
\begin{align*}
\left(\frac{1}{\varphi(Q)^{\theta}|Q|}\int_Q\omega(y)dy\right)^{\frac{1}{p}}
\left(\frac{1}{\varphi(Q)^{\theta}|Q|}\int_Q\omega(y)^{\frac{-1}{p-1}}dy\right)^
{\frac{1}{p'}}\leq C,
\end {align*}
for all cubes $Q$. Especially, when $p=1$, we obtain
\begin{align*}
\left(\frac{1}{\varphi(Q)^{\theta}|Q|}\int_Q\omega(y)dy\right)\leq 
C\mathop{\inf}\limits_{x\in {\rm Q}}\omega(x),
\end {align*}
where $Q(x,r)$ is the cube centered at $x$ with the sidelength $r$,  
$\varphi(Q):=(1+r)$, $|E|$ denotes the 
Lebesgue measure of $E$ and $\chi_E$ stands for the characteristic function of 
$E$. $f_Q$ denotes the average $f_Q:=\frac{1}{|Q|} \int_Qf(y)dy$, and the 
letter $ C $ will denote a constant not necessarily the same at each occurrence.
\begin{remark}
Let $A_p^{\infty}(\varphi):= {\textstyle \bigcup_{\theta \ge 0}} 
A_p^{\theta}(\varphi)$, and $A_{\infty}^{\infty}(\varphi):= {\textstyle 
\bigcup_{p\ge 1}} A_p^{\infty}(\varphi)$. The class $A_p^{\theta}(\varphi)$ is 
strictly larger than the class $A_p$ for all $1\leq p<\infty$. Moreover, taking  
for instance,  the weight $ \omega \in A_{p}^{\theta } (\varphi ) $, for $ 
\omega (x):=1+|x|^2 $, $ p> 1$, and $0< \theta \le 2n $, but $ \omega \notin 
A_{p} $. In particular, when $ \theta =0$, $A_p^0(\varphi)$ is equivalent 
to Muckenhoupt's class of weights $A_p$ for all $1\leq p<\infty$.
\end{remark}
In 2009, Lerner and Ombrosi et al \cite{b19} introduced a new class of weighted 
functions. The multiple weight $ A_{\vec p}$ is defined as follows.

Let $\vec{p}=(p_1,\ldots,p_m)$ and $1/p= {\textstyle \sum_{j=1}^{m}}  1/p_j$ 
with $1\leq p_1,\ldots,p_m<\infty$. Given 
$\vec{\omega}=(\omega_1,\ldots,\omega_m)$, each $\omega_j$ being nonnegative 
measurable, set $v_{\vec{\omega}}:=  {\textstyle 
\prod_{j=1}^{m}}\omega_j^{p/p_j}$. We say that $\vec\omega\in A_{\vec{p}}$ , if 
and only if
\begin{align*}
\mathop{\sup}\limits_{\rm 
Q}\left(\frac{1}{|Q|}\int_Qv_{\vec{\omega}}(x)dx\right)^{\frac1p}
\prod_{j=1}^{m}\left(\frac{1}{|Q|}\int_Q\omega_j(x)^{1-p'_j}dx\right)^{\frac{1}{p'_j}}<\infty,
\end {align*}
where the supremum is taken over all cubes $Q\subset \mathbb{R}^n$. When 
$p_j=1$, $j=1,\ldots, m$, the term $\left(\frac{1}{|Q|}\displaystyle \int_Q 
\omega_j(x)^{1-p'_j}dx\right)^{\frac{1}{p'_j}}$ is considered to be 
$\left(\mathop{\inf}\limits_{x\in {\rm Q}}\omega_j(x)\right)^{-1}$. 

In 2015, Bui \cite{b2} studied a new class of multiple weights $ A_{\vec 
p}^{\infty } (\varphi ) $. In 2015, Pan and Tang \cite{b25} obtained the 
pointwise estimates, strong type and weak end-point estmates for certain 
classes of multilinear operators and their iterated commutators with new $ BMO 
$ functions, where the weights are $ A_{\vec p}^{\infty } (\varphi ) $. For some other studies of $A_{\vec p}^{\infty}(\varphi)$, one can see, for 
instance, \cite{b38}.  

We say that $\vec\omega\in A_{\vec{p}}^{\theta}(\varphi)$, if and only if
\begin{align*}
\mathop{\sup}\limits_{\rm Q} 
\left(\frac{1}{\varphi(Q)^{\theta}|Q|}\int_Qv_{\vec{\omega}}(x)dx\right)^{\frac1p}
\prod_{j=1}^{m}\left(\frac{1}{\varphi(Q)^{\theta}|Q|}\int_Q\omega_j(x)^{1-p'_j}dx\right)^{\frac{1}{p'_j}}<\infty,
\nonumber
\end{align*}
where the supremum is taken over all cubes $Q\subset \mathbb{R}^n$ and  
$\theta\geq0$. When $p_j=1$, $j=1,\ldots, m$, the term 
$\left(\frac{1}{|Q|}\displaystyle \int_Q  
\omega_j(x)^{1-p'_j}dx\right)^{\frac{1}{p'_j}}$ is considered to be 
$\left(\mathop{\inf}\limits_{x\in {\rm Q}}\omega_j(x)\right)^{-1}$. 
\begin{remark}
For $1\leq p_1,\ldots,p_m<\infty$, set 
$A_{\vec{p}}^{\infty}(\varphi):= {\textstyle \bigcup_{\theta \ge 0}} 
A_{\vec{p}}^{\theta}(\varphi)$. When $\theta=0$, the class 
$A_{\vec{p}}^{0}(\varphi)$ is considered to be the class of multiple weights 
$A_{\vec{p}}$ introduced by {\rm \cite{b19}}.
\end{remark}

Inspired by the above multilinear square operator $ T $ and $A_{\vec p}^{\theta 
}(\varphi )$ weight functions, the new multilinear square operator with 
classical kernel is defined as follows. For any $ t\in (0,\infty ) $, let $ 
K_t(y_0,y_1,\ldots,y_m) $ be a locally integrable function defined away from 
the diagonal $y_0=y_1=\ldots=y_m$. It is called a classical kernel, if for some 
positive constants $ C $ and $\gamma  $, any $ N\ge 0 $,  it satisfies the 
following size condition,
\begin{equation}\label{1.2}	
\left (  \int_{0}^{\infty}|K_{t}(y_0, y_{1}, \ldots, y_{m})|^{2} \frac{d 
t}{t}\right ) ^{1 / 2} \leq 
\frac{C}{(\sum_{j=1}^m|y_0-y_j|)^{mn}(1+\sum_{j=1}^m|y_0-y_j|)^N},
\end{equation}
and the smoothness condition
\begin{align*}
&\left(\int_{0}^{\infty}\left|K_{t}\left(y_0, y_{1},\ldots, 
y_{m}\right)-K_{t}\left(y_0', y_{1},\ldots, y_{m}\right)\right|^{2} 
\frac{d t}{t}\right)^{1 / 2} \\
&\leq \frac{C|y_0-y_0'|^{\gamma }}{(\sum_{j=1}^m|y_0-y_j|)^{mn+\gamma 
}(1+\sum_{j=1}^m|y_0-y_j|)^N},
\end{align*}
whenever $|y_0-y_0'|\leq\frac{1}{2}\max_{1\leq j \leq m}|y_0-y_j|$.

Suppose that $\omega(t):[0,\infty) \mapsto  [0,\infty)$ is a nondecreasing 
function with $0<\omega(t)<\infty$. For $a>0$, one says that $\omega\in {\rm 
Dini}(a)$, if
\begin{align*}
[\omega]_{{\rm Dini}[a]}:=\int^{1}_{0}\frac{\omega^a(t)}{t}dt<\infty.  
\end{align*}
It is worth mentioning that Dini($a_1$) $\subset$ Dini($a_2$) when $0<a_1<a_2$.

The new multilinear square operator with kernel of Dini's type is defined as 
follows, where the size condition \eqref{1.2} remains unchanged and the 
smoothness condition is changed by
\begin{equation}\label{1.3}	
\begin{aligned}
&\left(\int_{0}^{\infty}\left|K_{t}\left(y_0, y_{1},\ldots, 
y_{m}\right)-K_{t}\left(y_0', y_{1}, \ldots, y_{m}\right)\right|^{2} 
\frac{d t}{t}\right)^{1 / 2} \\
&\leq 
\frac{C}{(\sum_{j=1}^m|y_0-y_j|)^{mn}(1+\sum_{j=1}^m|y_0-y_j|)^N}\omega 
\left( \frac{|y_0-y_0'|}{\sum_{j=1}^m |y_0-y_j|}\right),
\end{aligned}
\end{equation}
whenever $|y_0-y_0'|\leq\frac{1}{2}\max_{1\leq j \leq m}|y_0-y_j|$.

The conditions \eqref{1.3} can be further weakened, then a class of generalized 
kernels can be obtained as follows.
\begin{equation}\label{1.4}	
\begin{aligned}
&\left(\int_{0}^{\infty}\left(\int_{\left(\Delta_{k+2}\right)^{m} 
\backslash\left(\Delta_{k+1}\right)^{m}}\left|K_{t}\left(y_{0}, y_{1}, 
\ldots, y_{m}\right)-K_{t}\left(y_{0}^{\prime}, y_{1}, \ldots, 
y_{m}\right)\right|^{q} d \vec{y}\right)^{\frac{2}{q}} \frac{d 
t}{t}\right)^{\frac{1}{2}}\\
&\leq CC_k|y_0-y_o'|^{-\frac{mn}{q'}}\left(1+2^k |y_0-y_0'| 
\right)^{-N}2^{-\frac{kmn}{q'}},
\end{aligned}
\end{equation}
where $\Delta_k:=Q(y_0, 2^{k}\sqrt{mn}|y_0-y_0'|)$ denotes the cube centered at 
$y_0$ with the sidelength $ 2^{k}\sqrt{mn}|y_0-y_0'| $, $k\in \mathbb{N}$, $ 
1<  q\le 2$, $1/q+1/q'=1$ and  $C_k>0$. 


\begin{remark}
When $ \omega (t):=t^\gamma $ for some $ \gamma > 0 $ in the condition 
\eqref{1.3}, the new multilinear square operator with kernel of Dini's type is 
the new multilinear square operator with classical kernel.
\end{remark}


\begin{remark}
The class of generalized kernels includes the kernels of Dini's type, thus the 
new multilinear square operator with generalized kernel is larger than the new 
multilinear square operator with kernel of Dini's type. As a matter of fact, 
when $C_k:=\omega(2^{-k})$, it is easy to deduce that  the condition 
\eqref{1.3} 
implies the condition \eqref{1.4}, for any $1< q\le 2$.
\begin{align*}
&\left(\int_{0}^{\infty}\left(\int_{\left(\Delta_{k+2}\right)^{m} 
\backslash\left(\Delta_{k+1}\right)^{m}}\left|K_{t}\left(y_{0}, 
y_{1},\ldots, y_{m}\right)-K_{t}\left(y_{0}^{\prime}, y_{1}, \ldots, 
y_{m}\right)\right|^{q} d \vec{y}\right)^{\frac{2}{q}} \frac{d 
t}{t}\right)^{\frac{1}{2}}\\
&\leq \left[\int_{\left(\Delta_{k+2}\right)^{m} 
\backslash\left(\Delta_{k+1}\right)^{m}}\left(\int_{0}^{\infty}\left|K_{t}\left(y_{0},
y_{1},\ldots, y_{m}\right)-K_{t}\left(y_{0}^{\prime}, y_{1},\ldots, 
y_{m}\right)\right|^{q \cdot \frac{2}{q}} \frac{d 
t}{t}\right)^{\frac{q}{2}} d 
\vec{y}\right]^{\frac{1}{q}}\\
&\leq \Bigg(\int_{(\Delta_{k+2})^m \setminus 
(\Delta_{k+1})^m}\Bigg[\frac{C}{(\sum_{j=1}^m|y_0-y_j|)^{mn}(1+\sum_{j=1}^m|y_0-y_j|)^N}\\
&\ \ \ \times\omega \left( \frac{|y_0-y_0'|}{\sum_{j=1}^m 
|y_0-y_j|}\right)\Bigg]^q d\vec{y}\Bigg)^{1/q}\\
&\leq \Bigg(\int_{(\Delta_{k+2})^m \setminus 
(\Delta_{k+1})^m}\Bigg[\frac{C}{(\mathop{\max}\limits_{1\leq j \leq 
m}|y_0-y_j|)^{mn}(1+\mathop{\max}\limits_{1\leq j \leq 
m}|y_0-y_j|)^N}\\
&\ \ \ \times\omega \Bigg( 
\frac{|y_0-y_0'|}{\mathop{\max}\limits_{1\leq j \leq 
m}|y_0-y_j|}\Bigg)\Bigg]^q d\vec{y}\Bigg)^{1/q}\\
&\leq \left(\int_{(\Delta_{k+2})^m \setminus 
(\Delta_{k+1})^m}\left[\frac{C}{(2^k \sqrt{n}|y_0-y_0'|)^{mn}(1+2^k 
\sqrt{n}|y_0-y_0'|)^N}\omega \left( \frac{1}{2^k}\right)\right]^q 
d\vec{y}\right)^{1/q}\\
&\leq 
CC_k2^{\frac{-kmn}{q'}}|y_0-y_0'|^{\frac{-mn}{q'}}(1+2^k|y_0-y_0'|)^{-N}.
\end{align*}
\end{remark}

\begin{definition} \label{definition1.1}
Suppose that  $T$ is the multilinear square operator defined by \eqref{1.1}. If 
the following conditions are satisfied, $T$ is called the new multilinear 
square operator with generalized kernel.

\begin{itemize}
\item [\rm (i)] The kernel function satisfies the conditions \eqref{1.2} 
and \eqref{1.4} for any $N\geq 0$.
	
\item [\rm (ii)]  $T$ can be extended to be a 
bounded operator from $L^{s_1}\times \cdots \times L^{s_m}$ to $L^{s,\infty}$ 
for some $ 1\le s_j\le q', j=1,\ldots,m $ with	$1/s=1/s_1+\ldots+1/s_m$.
\end{itemize}
\end{definition}

\begin{remark}
The relationship between the new class of multilinear square operators and the  
classical multilinear square operators, and the relationship between new 
multiple weights and multiple weights are clearly explained in Figure 1.
\end{remark}

\begin{figure}[H]
\centering
\includegraphics[width=1\linewidth]{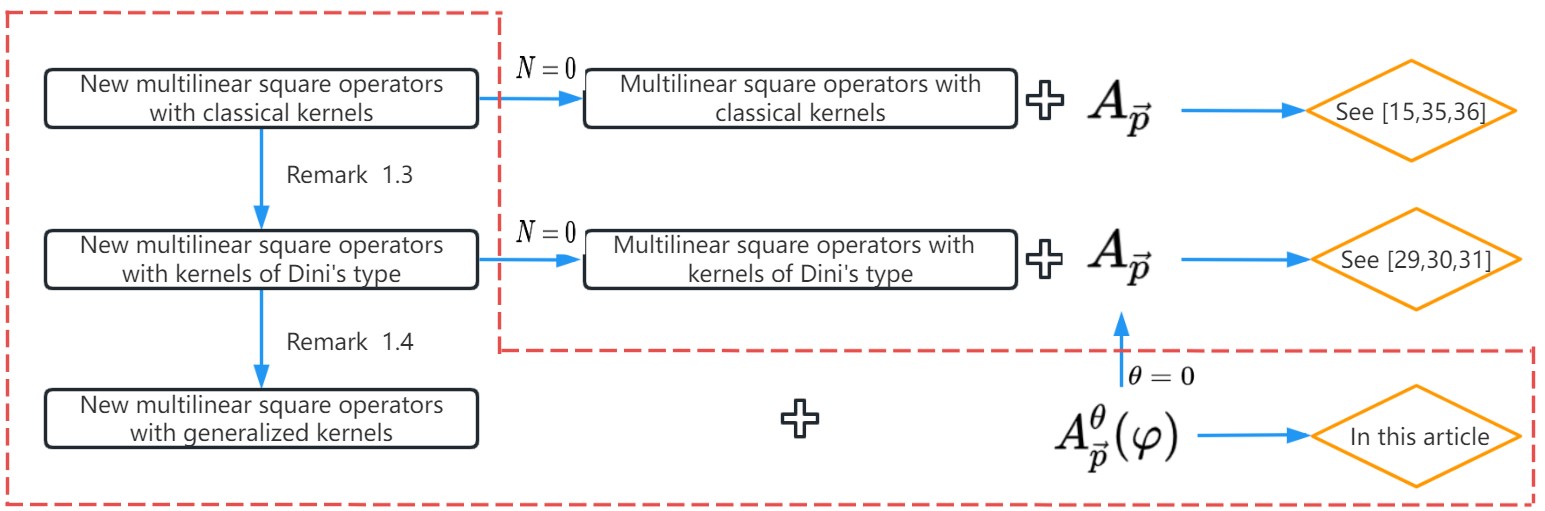} 
\caption{ Relation among of several classes of multilinear square operators and 
weight functions}
\label{FIG1}
\end{figure}

The structure of this article is as follows. In Section \ref{sec2}, the 
boundedness on the weighted Lebesgue spaces is obtained for the new multilinear 
square operator with generalized kernel. In Section \ref{sec3}, the boundedness 
on the weighted Morrey spaces is established, where the new Morrey space is 
acquired for the new multilinear square operator with generalized kernel.

\section{New weighted norm inequalities on Lebesgue spaces}\label{sec2}

\quad\quad The definition of the maximal function is introduced, as well as 
some properties of $ A_{p}^{\infty } (\varphi ) $ and $ A_{\vec p}^{\infty } 
(\varphi ) $ are given.

\subsection{Definitions and lemmas}

\begin{definition}
Let $0<\eta<\infty$, $ Q $ be a dyadic cube, and $f$ be a locally integral function. The dyadic maximal function 
$M_{\varphi,\eta}^\triangle$ is defined by 
\begin{align*}
	M_{\varphi,\eta}^\triangle f(x):=\mathop{\sup}\limits_{  {
			Q}\ni x}\frac{1}{\varphi(Q)^{\eta}|Q|}\int_Q|f(y)|dy.
\end{align*}

And the dyadic sharp maximal operator $M_{\varphi,\eta}^{\sharp,\triangle}$ is 
defined by 
\begin{align*}
	M_{\varphi,\eta}^{\sharp,\triangle} f(x)
	&:=\mathop{\sup}\limits_{{ 
			Q}\ni x,r<1}\frac{1}{|Q|}\int_{Q}|f(y)-f_Q|dy+\mathop{\sup}\limits_{{
			Q}\ni x,r\geq1}\frac{1}{\varphi(Q)^{\eta}|Q|}\int_{Q}|f(y)|dy\\
	&\simeq\mathop{\sup}\limits_{ { 
			Q}\ni x,r<1}\mathop{\inf}\limits_{ 
		c}\frac{1}{|Q|}\int_{Q}|f(y)-c|dy+\mathop{\sup}\limits_{{
			Q}\ni x,r\geq1}\frac{1}{\varphi(Q)^{\eta}|Q|}\int_{Q}|f(y)|dy.
\end{align*}
For $0<\delta < \infty$, denote as
\begin{align*}
	M^\triangle_{\delta,\varphi,\eta}f(x):=[M^\triangle_{\varphi,\eta}
	(|f|^\delta)]^{\frac1\delta}(x), M^{\sharp, 
		\triangle}_{\delta,\varphi,\eta}f(x):=[M^{\sharp, 
		\triangle}_{\varphi,\eta}(|f|^\delta)]^{\frac1\delta}(x).
\end{align*}
\end{definition}

\begin{lemma}\label{Lemma2.1}{\rm(\cite{b1})}   
The following statements are established.
\begin{itemize}
\item [\rm (i)] $A_p^{\infty}(\varphi) \subset A_q^{\infty}(\varphi) $, $1\leq 
p<q<\infty$.
\item [\rm (ii)] If $\omega \in A_p^{\infty}(\varphi )$ with $p>1$, then there 
exists a $\varepsilon>0$ such that $\omega \in  
A_{p-\varepsilon}^{\infty}(\varphi )$. Consequently, 
$A_p^{\infty}(\varphi ):= {\textstyle \bigcup_{q< p}} A_q^{\infty}(\varphi )$.
\item [\rm (iii)] If $\omega \in A_p^{\infty}(\varphi )$ for $p\geq1$, then 
there exist positive numbers $l$, $\delta$ and $C$ so that for all cubes 
$Q$,
\begin{align*}
\left(\frac{1}{|Q|}\int_Q\omega^{1+\delta}(x)dx\right)^{\frac{1}{1+\delta}}\leq
C\left(\frac{1}{|Q|}\int_Q\omega(x)dx\right)(1+r)^l.
\end{align*}
\end{itemize}
\end{lemma}

\begin{lemma}\label{Lemma2.2}{\rm(\cite{b32})}   The following statements are 
established.
\begin{itemize}
\item [\rm (i)] $\omega \in A_p^{\theta}(\varphi )$ if and only if 
$\omega^{1-p'} \in A_{p'}^{\theta}(\varphi )$,  
$1/p+1/p'=1$.
\item [\rm (ii)] If $\omega_0 , \omega_1 \in A_p^{\theta}(\varphi )$, $p\geq1$, 
then $\omega_0^{\alpha}\omega_1^{1-\alpha} \in A_p^{\theta}(\varphi )$ for 
any $0<\alpha<1$.
\item [\rm (iii)] If $\omega \in A_p^{\theta}(\varphi )$ for $1\leq p<\infty$, 
then
\begin{align*}
\frac{1}{\varphi(Q)^{\theta}|Q|}\int_Q |f(x)|dx \leq 
C\left(\frac{1}{\omega(5Q)}\int_Q|f(x)|^p\omega(x)dx\right)^{\frac1p}.
\end{align*}
\end{itemize}
Especially, let $f:=\chi_E$ for any measurable set $E\subset Q$,
\begin{align*}
\frac{|E|}{\varphi(Q)^{\theta}|Q|}\leq 
C\left(\frac{\omega(E)}{\omega(5Q)}\right)^{\frac{1}{p}}.
\end{align*}
\end{lemma}

\begin{lemma}\label{Lemma2.3}{\rm(\cite{b2})}  
Let $1\leq p_1,\ldots,p_m<\infty$ and 
$\vec{\omega}=(\omega_1,\ldots,\omega_m)$. So the following statements 
are equivalent.
\begin{itemize}
\item [\rm (i)] $\omega_j^{1-p'_j} \in A_{mp'_j}^{\infty},$ 
$j=1,\ldots,m $, and $v_{\vec{\omega}} \in A_{mp}^{\infty}(\varphi )$.
\item [\rm (ii)] $\vec\omega \in A_{\vec{p}}^{\infty}(\varphi ).$
\end{itemize}
\end {lemma}

\begin{remark}
The $A_{\vec{p}}^{\infty}(\varphi )$ weight function is not increasing which 
means for $\vec{p}=(p_1,\ldots,p_m)$ and $\vec{q}=(q_1,\ldots,q_m)$ with 
$p_j<q_j$, $j=1,\ldots, m$, the following may not be true 
$A_{\vec{p}}^{\infty}(\varphi ) \subset A_{\vec{q}}^{\infty}(\varphi )$.
\end{remark}

\begin{lemma}\label{Lemma2.4}{\rm(\cite{b32})}
Let $1<p<\infty$, $0<\eta<\infty$, $\omega\in A_p^\infty(\varphi )$ and 
$f\in L^p(\omega)$, the following statements hold.
\begin{align*}
\Vert f \Vert_{L^p(\omega)}\leq\Vert M_{\varphi,\eta}^\triangle f 
\Vert_{L^p(\omega)}\leq C \Vert M_{\varphi,\eta}^{\sharp, \triangle} f 
\Vert_{L^p(\omega)}.
\end{align*}
\end {lemma}

\begin{lemma} \label{Lemma2.5}{\rm(\cite{b32})}
Let $0<\eta<\infty$, $\omega\in A_\infty^\infty(\varphi )$ and 
$\delta>0$. Let $\psi: (0,\infty)\to(0,\infty)$ be doubling, which means that 
there is a constant $C_0$ and for $a>0$, $\psi(2a)\leq C_0\psi(a)$. Then there 
exists a constant C depending upon the doubling condition of  $\psi$ and the $ 
A_\infty^\infty(\varphi )$ condition of  $\omega$ such that,
\begin{align*}
\mathop{\sup}\limits_{\lambda>0}\psi(\lambda)\omega\left(\lbrace 
y\in\mathbb{R}^n, M_{\delta.\varphi,\eta}^\triangle f(y)>\lambda\rbrace\right)
\leq C \mathop{\sup}\limits_{\lambda>0}\psi(\lambda)\omega\left(\lbrace 
y\in\mathbb{R}^n, M_{\delta.\varphi,\eta}^{\sharp, \triangle} 
f(y)>\lambda\rbrace\right).
\end{align*}
\end {lemma}

For $0<\eta<\infty$, $0<\delta < \infty$, $\vec{f}=(f_1,\ldots, f_m)$, the 
multilinear maximal operators $\mathcal{M}_{\varphi,\eta}$ and $ 
\mathcal{M}_{\delta,\varphi,\eta} $ are defined by
\begin{align*}
	\mathcal{M}_{\varphi,\eta}(\vec{f})(x):=\mathop{\sup}\limits_{ {
			Q}\ni x}\prod_{j=1}^m\frac1{\varphi(Q)^\eta \left | Q \right | } \int_Q 
	|f_j(y_j)| dy_j,
\end{align*}
and
\begin{align*}
	\mathcal{M}_{\delta,\varphi,\eta}(\vec f)(x)
	:=\left [ \mathcal{M}_{\varphi,\eta}(|\vec{f}|^{\delta })(x) \right ] 
	^{\frac{1}{\delta } }
	=\left(\mathop{\sup}\limits_{ 
		{ Q}\ni x}\prod_{j=1}^m\frac1{\varphi(Q)^\eta \left | Q \right |} \int_Q 
	|f_j(y_j)|^\delta dy_j \right)^{\frac1\delta}.
\end{align*}

\begin{lemma}\label{Lemma2.6}{\rm(\cite{b25})}
Let $\vec \omega\in A_{\vec{p}}^\infty(\varphi )$ and $1<p_j<\infty$, 
$j=1,\ldots, m$, $1/p = {\textstyle \sum_{j=1}^{m}}  1/p_j
$, then there exists some $\eta_0>0$ depending on $p, m , p_j$, such that
\begin{align*}
\Vert \mathcal{M}_{\varphi,\eta_0}(\vec{f}) \Vert_{L^{p}(v_{\vec \omega })}\leq 
C\prod_{j=1}^{m}\Vert f_j \Vert_{L^{p_j}(\omega_j)}.
\end{align*}
\end {lemma}

\begin{lemma}\label{Lemma2.7}{\rm(\cite{b25})}
Let $\vec \omega\in A_{\vec p}^{\infty } (\varphi )$ and $1\leq p_j<\infty$, 
$1/p = {\textstyle \sum_{j=1}^{m}}  1/p_j$, 
$j=1,\ldots, m$, then there exists some $\theta_0>0$ depending on $p, m, p_j$, 
such that
\begin{align*}
\Vert \mathcal{M}_{\varphi,\theta_0}(\vec{f}) \Vert_{L^{p,\infty}(v_{\vec 
\omega })}\leq C\prod_{j=1}^{m}\Vert f_j \Vert_{L^{p_j}(\omega_j)}.
\end{align*}
\end {lemma}

\begin{lemma}\label{Lemma2.8}{\rm(\cite{b19})}
Suppose $0<p<q<\infty$. Then there is a 
positive constant $C:=C_{p,q}$, such that for all measurable functions $f$,
\begin{align*}
|Q|^{-1/p}\Vert f \Vert_{L^p (Q)}\leq C|Q|^{-1/q}\Vert f 
\Vert_{L^{q,\infty}(Q)}.
\end{align*}
\end {lemma}


\subsection{Main theorem and its proof}

\quad\quad The boundedness on the weighted Lebesgue spaces are obtained for new 
multilinear square operators with generalized kernels. The main theorems of this 
section is given below.

\begin{theorem}\label{Theorem2.1} 
Suppose that $m\geq2$, $T$ is the new multilinear square operator with 
generalized kernel as in Definition \ref{definition1.1}, $ {\textstyle 
\sum_{k=1}^{\infty }}  C_{k} < \infty $, $ \vec p=(p_1, \ldots ,p_m) $, $\vec 
\omega=(\omega_1,\ldots ,\omega_m) \in 
A_{\vec{p}/q'}^\infty(\varphi )$, $v_{\vec\omega}= {\textstyle \prod_{j=1}^{m}} 
\omega_j^{p/p_j}$ and $1/p = {\textstyle \sum_{j=1}^{m}}  1/p_j$. Then the 
following statements hold.

\begin{itemize}
\item [\rm (i)] If $q'<p_j<\infty, j=1,\ldots,m$, then there exists a 
constant $C>0$ such that
\begin{align*}
\Vert T(\vec{f}) \Vert_{L^{p}(v_{\vec\omega})}\leq 
C\prod_{j=1}^{m}\Vert f_j \Vert_{L^{p_j}(\omega_j)}.
\end{align*}

\item [\rm (ii)] If $q'\leq p_j<\infty, j=1,\ldots,m$, and at least one of 
$p_j=q'$, then there exists a constant $C>0$ such that
\begin{align*}
\Vert T(\vec{f}) \Vert_{L^{p,\infty}(v_{\vec\omega})}\leq 
C\prod_{j=1}^{m}\Vert f_j \Vert_{L^{p_j}(\omega_j)}.
\end{align*}
\end{itemize}
\end{theorem}

\noindent{\it \rm Proof.} 
Firstly, we will establish the point estimate of  the sharp maximal function.  
For $0<\delta<1/m$ and $0<\eta<\infty$, there exists a constant $C>0$ 
such that for all bounded measurable functions $\vec 
f=(f_1,\ldots,f_m)$ with compact support,
\begin{equation}\label{2.1}
M_{\delta,\varphi,\eta}^{\sharp, \triangle}(T(\vec f))(x)\leq 
C\mathcal{M}_{q',\varphi,\eta}(\vec f)(x).
\end{equation}

Fixed $x\in\mathbb{R}^n$, for any dyadic cube $ Q:=Q(x_0,r)\ni x$, to prove the 
formula \eqref{2.1}, we consider two cases about $ r $: $r<1$ and $r\geq1$.

\noindent{\it \rm \textbf{Case 1}:} $r<1$. Since $0<\delta<1/m< 1$
and $\left | \left | a_{1}  \right | ^{t} - | \left | a_{2}  \right | ^{t} 
\right | \le \left | a_{1}- a_{2}  \right | ^{t} $ for $0<t<1$, for any 
constant $A$, we have
\begin{align*}
\left(\frac1{|Q|}\int_Q \left||T(\vec f)(z)|^\delta-|A|^\delta \right| 
dz\right)^{\frac1\delta}\leq
\left(\frac1{|Q|}\int_Q \left|T(\vec f)(z)-A \right|^\delta 
dz\right)^{\frac1\delta}.
\end{align*}

Let $Q^*:=14n\sqrt{nm}Q$, we split each $f_j$ as 
$f_j=f_j^0+f_j^\infty$, for each $f_j^0=f_j \chi_{Q^*}$, $ f_j^\infty=f_j-f_j^0 
$. This yields
\begin {align*}
\prod_{j=1}^mf_j(y_j)=\prod_{j=1}^m f_j^{0}(y_j)+\sum_{(\alpha_1,\ldots,\alpha_m)\in 
\mathscr{L}}f_1^{\alpha_1}(y_1)\cdots f_m^{\alpha_m}(y_m),
\end {align*}
where  $\mathscr{L}:=\{(\alpha_1,\ldots,\alpha_m)$: there is at least one 
$\alpha_j=\infty\}$.

Taking $z_0\in 4Q\setminus 3Q$ and $A:=\left(\int_{0}^{\infty} \left | \int_{\left(\mathbb{R}^{n}\right)^{m} 
}K_{t}\left(z_{0} 
,\vec{y}\right) \sum_{(\alpha_1,\ldots,\alpha_m)\in \mathscr{L}}\prod_{j=1}^{m} 
f_{j}^{\alpha_j}\left(y_{j}\right) d \vec{y} \right | ^{2} \frac{d 
	t}{t}\right)^{\frac{1}{2}}$, we conclude that for any $z\in Q$,
	
\begin  {align*}
& |T(\vec f)(z)-A| \\
 &\leq  \Bigg(\int_{0}^{\infty} \bigg | 
 \int_{\left(\mathbb{R}^{n}\right)^{m} }\Big(K_{t}(z, \vec{y})\big ( \prod_{j=1}^m f_j^{0}(y_j)+\sum_{(\alpha_1,\ldots,\alpha_m)\in 
 	\mathscr{L}}f_1^{\alpha_1}(y_1)\cdots f_m^{\alpha_m}(y_m) \big ) \\
&\ \ \ \ -K_{t}\left(z_{0} 
,\vec{y}\right) \sum_{(\alpha_1,\ldots,\alpha_m)\in \mathscr{L}}\prod_{j=1}^{m} 
f_{j}^{\alpha_j}\left(y_{j}\right)\Big)  d \vec{y} \bigg | ^{2} \frac{d 
	t}{t}\Bigg)^{\frac{1}{2}} \\
&\leq T(f_1^0,\ldots,f_m^0)(z)\\
&\ \ \ \ + \sum_{(\alpha_1,\ldots,\alpha_m)\in 
	\mathscr{L}}\left(\int_{0}^{\infty} \left | 
\int_{\left(\mathbb{R}^{n}\right)^{m} 
}\left(K_{t}(z, \vec{y})-K_{t}\left(z_{0} 
,\vec{y}\right)\right) \prod_{j=1}^{m} 
f_{j}^{\alpha_j}\left(y_{j}\right) d \vec{y} \right | ^{2} \frac{d 
	t}{t}\right)^{\frac{1}{2}}.
\end {align*}
We have
\begin  {align*}
&\left(\frac1{|Q|}\int_Q |T(\vec f)(z)-A|^\delta dz\right)^{\frac1\delta} \\
&\leq C \left(\frac1{|Q|}\int_Q |T(f_1^0,\ldots,f_m^0)(z)|^\delta 
dz\right)^{\frac1\delta}\\
&\ \ \ \ + C \sum_{(\alpha_1,\ldots,\alpha_m)\in \mathscr{L}}
\frac{1}{|Q|} \int_{Q}\left(\int_{0}^{\infty}\Bigg| 
\int_{\left(\mathbb{R}^{n}\right)^{m} 
}(K_{t}(z, \vec{y})-K_{t}\left(z_{0} 
,\vec{y}\right))  \prod_{j=1}^{m} 
f_{j}^{\alpha _j}\left(y_{j}\right) d \vec{y} \Bigg| ^{2} \frac{d 
	t}{t}\right)^{\frac{1}{2}} d z\\
&:=\uppercase\expandafter{\romannumeral+1}+\uppercase\expandafter{\romannumeral+2}\\
&:=\uppercase\expandafter{\romannumeral+1}+\sum_{(\alpha_1,\ldots,\alpha_m)\in \mathscr{L}}\uppercase\expandafter{\romannumeral+2}_{\alpha _1\dots \alpha _m}.
\end {align*}

For $1\le s_1,\dots ,s_m\le q'$, as $T:L^{s_1}\times\cdots\times L^{s_m}\to L^{s,\infty}$, using Lemma 
\ref{Lemma2.8} and H\"{o}lder's inequality, we deduce that
\begin {align*}
\uppercase\expandafter{\romannumeral+1}
&\leq C\Vert T(f_1^0,\ldots,f_m^0) \Vert_{L^{s,\infty}(Q,\frac{dx}{|Q|})}
\leq 
C\prod_{j=1}^m\left(\frac1{|Q^*|}\int_{Q^*}|f_j(y_j)|^{s_j}dy_j\right)^{\frac1{s_j}}\\
&\leq 
C\prod_{j=1}^m\left(\frac1{|Q^*|}\int_{Q^*}|f_j(y_j)|^{q'}dy_j\right)^{\frac1{q'}}
\leq 
C\prod_{j=1}^m\left(\frac1{\varphi(Q^*)^\eta|Q^*|}\int_{Q^*}|f_j(y_j)|^{q'}dy_j\right)^{\frac1{q'}}\\
&\leq C\mathcal{M}_{q',\varphi,\eta}(\vec f)(x).
\end {align*}

To estimate the term $\uppercase\expandafter{\romannumeral+2}$, let    
  $\Delta_k:=Q(z_0, 2^{k}\sqrt{mn}|z-z_0|)$, $  
z\in Q $, $k\in \mathbb{N}_+$. Since  $\Delta_2 \subset Q^*$, $z_0\in 4Q\setminus 3Q$, $  
z\in Q $, we have $(\mathbb{R}^n)^m\setminus 
(Q^*)^m\subset (\mathbb{R}^n)^m\setminus (\Delta_2)^m$, $|z-z_0|\sim r$. Because each term of 
$\mathscr{L}$ contains at least one $\alpha_j=\infty, j=1, \ldots,m$, we choose 
one of them for discussion. Taking $N\geq m\eta/q'$, it follows that
\begin{align*}
\uppercase\expandafter{\romannumeral+2}_{\alpha _1\dots \alpha _m}
&\leq\frac{1}{|Q|} \int_{Q}\left(\int_{0}^{\infty} \left ( 
\int_{\left(\mathbb{R}^{n}\right)^{m} 
	\backslash\left(Q^{*}\right)^{m}}\left|K_{t}(z, \vec{y})-K_{t}\left(z_{0} 
,\vec{y}\right)\right|  \prod_{j=1}^{m} 
\left|f_{j}\left(y_{j}\right)\right| d \vec{y} \right ) ^{2} \frac{d 
	t}{t}\right)^{\frac{1}{2}} d z \\
&\leq \frac{1}{|Q|} \int_{Q} \sum_{k=1}^{\infty }\Bigg(\int_{0}^{\infty}  \Big
( \int_{\left(\Delta_{k+2}\right)^{m} 
	\backslash\left(\Delta_{k+1}\right)^{m}}\left|K_{t}(z, 
\vec{y})-K_{t}\left(z_{0} ,\vec{y}\right)\right| \\
&\ \ \ \ \times \prod_{j=1}^{m} 
\left|f_{j}\left(y_{j}\right)\right| d \vec{y} \Big) ^{2}  \frac{d 
	t}{t}\Bigg)^{\frac{1}{2}} d z\\
&\leq \frac{1}{|Q|} \int_{Q} \sum_{k=1}^{\infty } \Bigg(  \int_{0}^{\infty}     
\Big( \int_{\left(\Delta_{k+2}\right)^{m} 
\backslash\left(\Delta_{k+1}\right)^{m}}\left|K_{t}(z, 
\vec{y})-K_{t}\left(z_{0} ,\vec{y}\right)\right|^{q}  d \vec{y}  \Big) 
^\frac{2}{q}\\
&\ \ \ \ \times  \Big( \int_{\left(\Delta_{k+2}\right)^{m} 
} \prod_{j=1}^{m} 
\left|f_{j}\left(y_{j}\right)\right|^{q'}  d \vec{y} \Big) ^\frac{2}{q'}   
\frac{d t}{t}\Bigg) ^{\frac{1}{2}} d z\\
&\leq \frac{1}{|Q|} \int_{Q} \sum_{k=1}^{\infty } \Bigg(  \int_{0}^{\infty}     
\Big( \int_{\left(\Delta_{k+2}\right)^{m} 
\backslash\left(\Delta_{k+1}\right)^{m}}\left|K_{t}(z, 
\vec{y})-K_{t}\left(z_{0} ,\vec{y}\right)\right|^{q}  d \vec{y}  \Big) 
^\frac{2}{q} \frac{d t}{t}\Bigg) ^{\frac{1}{2}} \\
&\ \ \ \ \times  \prod_{j=1}^{m}\Big( \int_{\Delta_{k+2}}  
\left|f_{j}\left(y_{j}\right)\right|^{q'}  d y_j \Big) ^\frac{1}{q'}dz\\
&\leq \frac{C}{|Q|}\int_Q \sum_{k=1}^\infty C_k 
2^{\frac{-kmn}{q'}}|z-z_0|^{\frac{-mn}{q'}}(1+2^k|z-z_0|)^{-N}  
\prod_{j=1}^m\left(\int_{\Delta_{k+2}}|f_j(y_j)|^{q'} dy_j 
\right)^{\frac 1{q'}} dz\\
&\leq \frac{C}{|Q|}\int_Q \sum_{k=1}^\infty \frac{C_k 
\varphi(\Delta_{k+2})^{m\eta/q'}|\Delta_{k+2}|^{m/q'}}{(1+2^kr)^N(2^kr)^{mn/q'}}\\
&\ \ \ \ \times \prod_{j=1}^m\left(\frac{1}{\varphi(\Delta_{k+2})^{\eta}|\Delta_{k+2}|}
\int_{\Delta_{k+2}}|f_j(y_j)|^{q'}
dy_j\right)^{\frac 1{q'}}dz\\
&\leq \frac{C}{\left | Q \right | } \int_Q \sum_{k=1}^\infty C_{k}
(1+2^{k}r)^{\frac{m\eta }{q'}- N }  \mathcal{M}_{q',\varphi,\eta}(\vec 
f)(x) dz\\
&\leq C\mathcal{M}_{q',\varphi,\eta}(\vec f)(x).
\end{align*}

Thus, $\uppercase\expandafter{\romannumeral+2}\leq 
C\mathcal{M}_{q',\varphi,\eta}(\vec f)(x).$

\noindent{\it \rm \textbf{Case 2}:} $r\geq1$. Let $\widetilde{Q^*}=8Q$, we split each $f_j$ as $f_j=f_j^0+f_j^\infty$, for each $f_j^0=f_j \chi_{\widetilde{Q^*}}$, 
$f_j^\infty=f_j-f_j^0 $. There is
\begin {align*}
\prod_{j=1}^mf_j(y_j)
&=\prod_{j=1}^m f_j^{0}(y_j)+\sum_{(\alpha_1,\ldots,\alpha_m)\in 
	\mathscr{L}}f_1^{\alpha_1}(y_1)\cdots f_m^{\alpha_m}(y_m),
\end {align*}
where  $\mathscr{L}:=\{(\alpha_1,\ldots,\alpha_m)$: there is at least one 
$\alpha_j=\infty\}$.

We have
\begin {align*}
\left(\frac1{\varphi(Q)^\eta|Q|}\int_Q |T(\vec f)(z)|^\delta 
dz\right)^{\frac1\delta}
&\leq \frac{C}{\varphi(Q)^{\eta/\delta}}\left(\frac1{|Q|}\int_Q 
|T(f_1^0,\ldots,f_m^0)(z)|^\delta dz\right)^{\frac1\delta}\\
&\ \ \ \ +C\sum_{\alpha_1\ldots\alpha_m\in 
	\mathscr{L}}\frac{1}{\varphi(Q)^{\eta/\delta}}\left(\frac1{|Q|}\int_Q|T(f_1^{\alpha_1},\ldots,f_m^{\alpha_m})(z)|^\delta
 dz\right)^{\frac1\delta}\\
&:=\uppercase\expandafter{\romannumeral+3}+\uppercase\expandafter{\romannumeral+4}\\
&:=\uppercase\expandafter{\romannumeral+3}+\sum_{\alpha_1\ldots\alpha_m\in 
	\mathscr{L}}\uppercase\expandafter{\romannumeral+4}_{\alpha _1\dots \alpha _m}.
\end {align*}

To estimate $\uppercase\expandafter{\romannumeral+3}$. For $1\le s_1,\dots ,s_m\le q'$, according to $T:L^{s_1}\times\cdots\times L^{s_m}\to L^{s,\infty}$, using Lemma 
\ref{Lemma2.8} and H\"{o}lder's inequality, we deduce that
{\setlength\abovedisplayskip{10pt}
\setlength\belowdisplayskip{10pt}
\begin {align*}
\uppercase\expandafter{\romannumeral+3}
&\leq C\frac{1}{\varphi(Q)^{m\eta}}\Vert T(f_1^0,\ldots,f_m^0) 
\Vert_{L^{s,\infty}(Q,\frac{dx}{|Q|})}
\leq 
C\frac{1}{\varphi(Q)^{m\eta}}\prod_{j=1}^m\left(\frac1{|Q^*|}\int_{Q^*}|f_j(y_j)|^{s_j}dy_j\right)^{\frac1{s_j}}\\
&\leq 
C\frac{1}{\varphi(Q)^{m\eta}}\prod_{j=1}^m\left(\frac1{|Q^*|}\int_{Q^*}|f_j(y_j)
|^{q'}dy_j\right)^{\frac1{q'}}
\leq 
C\prod_{j=1}^m\left(\frac1{\varphi(Q^*)^\eta|Q^*|}\int_{Q^*}|f_j(y_j)|^{q'}dy_j\right)^{\frac1{q'}}\\
&\leq C\mathcal{M}_{q',\varphi,\eta}(\vec f)(x).
\end {align*}

For the term $\uppercase\expandafter{\romannumeral+4}$, when $(y_1,\ldots,y_m)\in (2^{k+3}Q)^m \setminus 
(2^{k+2}Q)^m$, $z\in Q, z_0\in 4Q\setminus 3Q$, we obtain $ {\textstyle \sum_{j=1}^{m}} 
|z-y_j|\sim 2^{k}r$, $|z-z_0|\sim r$. Taking $N >m\eta/q'$, by Minkowski's inequality 
and the size condition \eqref{1.2}, it follows that
{\setlength\abovedisplayskip{10pt}
\setlength\belowdisplayskip{10pt}
\begin {align*}
\uppercase\expandafter{\romannumeral+4}_{\alpha_1\ldots\alpha_m}
&\leq \frac{C}{\varphi (Q)^{\eta /\delta }|Q|} \int_{Q}\left(\int_{0}^{\infty}  \bigg( 
\int_{\left(\mathbb{R}^{n}\right)^{m} 
	\backslash(\widetilde{Q^*})^{m}}K_{t}(z, \vec{y})
\prod_{j=1}^{m} 
f_{j}\left(y_{j}\right) d \vec{y}  \bigg) ^{2} \frac{d 
	t}{t}\right)^{\frac{1}{2}} d z\\
&\leq \frac{C}{\varphi (Q)^{\eta /\delta }|Q|} \int_{Q} \int_{\left(\mathbb{R}^{n}\right)^{m} 
	\backslash(\widetilde{Q^*})^{m}}\left(\int_{0}^{\infty} \left |  
K_{t}(z, \vec{y})\right|^{2} \frac{d 	t}{t}\right)^{\frac{1}{2} }
\prod_{j=1}^{m} \left|f_{j}\left(y_{j}\right)\right| d \vec{y} d z\\
&\leq \frac{C}{\varphi (Q)^{\eta /\delta }|Q|} \int_{Q} \sum_{k=1}^{\infty } 
\int_{(2^{k+3}Q)^m\setminus (2^{k+2}Q)^m}\frac{\prod_{j=1}^{m} \left|f_{j}\left(y_{j}\right)\right| d 
	\vec{y} dz}{\left ( 
	\sum_{j=1}^{m} \left | z- y_{j}  \right |  \right ) ^{mn}\left ( 1+  
	\sum_{j=1}^{m} \left | z- y_{j}  \right |  \right ) ^{N} }  \\
&\leq  \frac{C}{\varphi (Q)^{\eta /\delta }|Q|}\int_Q \sum_{k=1}^\infty 
\frac{|2^{k+3}Q|^m}{(2^kr)^{mn}(1+2^kr)^{N}}\prod_{j=1}^m\left(\frac{1}
{|2^{k+3}Q|}\int_{2^{k+3}Q}|f_j(y_j)|dy_j\right)dz\\
&\leq  \frac{C}{\varphi (Q)^{\eta /\delta }|Q|}\int_Q \sum_{k=1}^\infty 
\frac{|2^{k+3}Q|^{m}\varphi(2^{k+3}Q)^{m\eta/q'}}{(2^kr)^{mn}
	(1+2^kr)^N}\\
&\ \ \ \ \times \prod_{j=1}^m\left(\frac{1}{\varphi(2^{k+3}Q)^{\eta}|2^{k+3}Q|}
\int_{2^{k+3}Q}|f_j(y_j)|^{q'}dy_j\right)^{\frac 1{q'}} dz\\
&\leq  \frac C{|Q|}\int_Q \sum_{k=1}^\infty (2^k)^{\frac{m\eta }{q'} - N} 
\prod_{j=1}^m\left(\frac{1}{\varphi(2^{k+3}Q)^{\eta}|2^{k+3}Q|}
\int_{2^{k+3}Q}|f_j(y_j)|^{q'}dy_j\right)^{\frac 1{q'}} dz\\
&\leq C\mathcal{M}_{q',\varphi,\eta}(\vec f)(x).
\end{align*}

Thus, $\uppercase\expandafter{\romannumeral+4}\leq 
C\mathcal{M}_{q',\varphi,\eta}(\vec f)(x)$.
By the above proof, the formula \eqref{2.1} was obtained. 

Secondly, we prove the strong type boundedness of the new multilinear square 
operator with generalized kernel on the weighted Lebesgue spaces. 

Choose a $ \delta  $ such that $ 0< \delta < 1/m $, and $ \eta =\eta _0 $ in Lemma 
\ref{Lemma2.6} for $ \vec \omega \in A_{\vec p/q'}^\infty(\varphi ) $, $p_{j}/q'> 1$. It 
derives from Lemma \ref{Lemma2.3} that $v_{\vec\omega}\in 
A_{mp/q'}^\infty(\varphi ) \subset A_{p/\delta } ^\infty(\varphi )$. Using 
Lemma \ref{Lemma2.4}, the formula \eqref{2.1} and Lemma \ref{Lemma2.6}, we 
deduce that
\begin{align*}
\| T(\vec f) \|_{L^p(v_{\vec\omega})} 
&=\left \| |T(\vec f)|^{\delta } \right \| _{L^{p/\delta  
}(v_{\vec\omega})}^{1/\delta  } 
\leq \| M_{\varphi,\eta _0}^{\bigtriangleup }(|T(\vec f)|^{\delta }) 
\|_{L^{p/\delta }(v_{\vec\omega})}^{1/\delta }\\
&\leq C \| M_{\delta,\varphi,\eta _0}^{\sharp, \triangle}(T(\vec f)) 
\|_{L^p(v_{\vec\omega})}
\leq C\| \mathcal{M}_{q',\varphi,\eta _0}(\vec f) \|_{L^p(v_{\vec\omega})}\\
&=C\| \mathcal{M}_{\varphi,\eta _0}(|\vec f|^{q'}) 
\|_{L^{p/q'}(v_{\vec\omega})}^{1/q'}
\leq C\prod_{j=1}^{m}\| |f_j|^{q'} 
\|_{L^{p_j/q'}(\omega_j)}^{1/q'}\\
&=C \prod_{j=1}^{m}\| f_j \|_{L^{p_j}(\omega_j)}.
\end{align*}

The expected result is obtained.

Finally, we establish the weak type estimate for the new multilinear square
operator $ T $ with generalized kernel on the weighted Lebesgue spaces.
\begin{align*}
 \| T( \vec f )   \|_{L^{p,\infty } \left ( \upsilon 
_{\vec \omega }  \right ) }
&=\sup_{\lambda > 0}\lambda \left | \nu _{\vec\omega }\left \{ x\in 
\mathbb{R}^n: \left | T ( \vec f  ) ( x  )  \right 
|>\lambda  \right \} 
\right |^\frac{1}{p} \\
&=\sup_{\lambda > 0}\left [\lambda^\delta  \left |\nu _{\vec\omega }\left \{ 
x\in \mathbb{R}^n: \left |T ( \vec f  ) ( x  ) \right 
|^\delta >\lambda ^\delta  \right \} \right |^\frac{\delta }{p}\right 
]^\frac{1}{\delta } \\
&=\sup_{\lambda > 0}\left [\lambda \left |\nu _{\vec\omega }\left \{ x\in 
\mathbb{R}^n: \left | T ( \vec f  ) ( x  ) \right 
|^\delta >\lambda   \right \} \right |^\frac{\delta }{p}\right 
]^\frac{1}{\delta } \\
&=\left \|T( \vec f  ) ^\delta  \right \|_{L^{\frac{p}{\delta } 
		,\infty } \left ( \upsilon _{\vec \omega }  \right ) }^{\frac{1}{\delta } } 
\le \left \| M_{\varphi ,\theta _0}^{\bigtriangleup} \left ( \left | 
T( \vec f  ) \right | ^\delta \right )   \right \|_{L^{\frac{p}{\delta } 
		,\infty } \left ( \upsilon _{\vec 
		\omega }  \right ) }^{\frac{1}{\delta } } \\
&=\sup_{\lambda > 0}\left [\lambda \left |\nu _{\vec\omega }\left \{ x\in 
\mathbb{R}^n: M_{\varphi ,\theta _0}^{\bigtriangleup} \left ( \left 
|  T ( \vec f  ) ( x  ) \right |^\delta \right )  >\lambda   \right \} \right 
|^\frac{\delta }{p}\right ]^\frac{1}{\delta }\\
&=\sup_{\lambda > 0}\lambda \left |\nu _{\vec\omega }\left 
\{ x\in \mathbb{R}^n: M_{\varphi ,\theta _0}^{\bigtriangleup}  \left 
( \left |T ( \vec f  ) ( x  ) \right |^\delta  \right ) >\lambda ^\delta    
\right \} \right |^\frac{1 }{p} \\
&=\sup_{\lambda > 0}\lambda \left |\nu _{\vec\omega }\left 
\{ x\in \mathbb{R}^n: \left (M_{\varphi ,\theta _0}^{\bigtriangleup}  
\left ( \left |T ( \vec f  ) ( x  ) \right |^\delta \right )  \right 
)^\frac{1}{\delta }  >\lambda     \right \} \right |^\frac{1 }{p} \\
&:=H ,
\nonumber
\end{align*}
where $ \theta _0 $ is taken by Lemma \ref{Lemma2.7} for $ \vec \omega \in 
A_{\vec p/q'}^\infty(\varphi ) $.

Next, we estimate $ H  $. Let $ \psi   (\lambda ):=\lambda ^{p} $, then $ 
\psi   (2\lambda )=2^{p} \lambda ^{p} \le 2^{p}\psi  (\lambda ) $, so $ 
\psi  $ satisfies the double condition. Using Lemma \ref{Lemma2.5}, Lemma 
\ref{Lemma2.7} and \eqref{2.1}, we have
\begin{align*}
H
&\le C \sup_{\lambda > 0}\lambda \left |\nu _{\vec\omega }\left 
\{ x\in \mathbb{R}^n:M_{\delta,\varphi,\theta 
_0}^{\sharp,\bigtriangleup}  
\left ( T ( \vec f  ) \right )  ( x  ) >\lambda   \right \} \right |^\frac{1 
}{p} \\
&=C\left \| M_{\delta ,\varphi ,\theta _0}^{\sharp ,\bigtriangleup} 
\left ( T( \vec f  ) \right )  \right \|_{L^{p,\infty } \left ( \upsilon _{\vec \omega }  \right 
	) }\le C\left \| \mathcal{M}_{q',\varphi ,\theta _0} ( \vec f  ) 
\right \|_{L^{p,\infty } \left ( \upsilon _{\vec \omega }  \right ) } \\
&=C \left \|\mathcal{M}_{\varphi ,\theta _0}( |\vec f|  ^{q'}) \right 
\|_{L^{p/q' ,\infty } \left ( \upsilon _{\vec \omega }  \right ) 
}^{1/q' } \le C\prod_{j=1}^{m} \left \||  f_j  | ^{q'} \right 
\|_{L^{p_j/q' } ( \omega  _{j }   ) }^{1/q' }  \\
&=C\prod_{j=1}^{m} \left \|  f_j   \right \|_{L^{p_j} ( \omega  _{j }   ) }.
\nonumber
\end{align*}

The proof of Theorem \ref{Theorem2.1} is finished.

\section{New weighted norm inequalities on new Morrey spaces}\label{sec3}

\quad\quad In 1938, Morrey \cite{b22} first introduced the classical Morrey 
spaces to investigate the local behavior of solutions to second-order elliptic 
partial 
differential equations. In recent years, many authors have studied the 
boundedness of operators on Morrey spaces and many properties of solutions to 
PDEs are concerned with the boundedness of operators on Morrey spaces. In 
2015, Pan and Tang \cite{b24} studied the boundedness for some Schr\"{o}dinger 
type operators on weighted Morrey spaces. In 2018, Trong and Truong \cite{b34} 
established  the weighted strong type and weak type estimates for the Riesz 
transforms and fractional integrals associated to Schr\"{o}dinger operators. 
For some other studies of weighted Morrey spaces, one can refer to \cite{b18,b26}.
Next, the definition of weighted Morrey spaces is given as follows.
\begin{definition}{\rm(\cite{b34})}
Let $ u $, $ w $ be two weights and $ \lambda \in [0,1) $, $ 1\le p< l\le 
\infty  $, $ \alpha \in (- \infty ,\infty ) $, $ Q:=Q(z,r)$. The strong Morrey 
space $ M_{\alpha ,\lambda }^{p,l} (u,w) $ is defined as the set of all 
measurable functions $ f $ satisfying $ \left \| f \right \| _{M_{\alpha 
,\lambda }^{p,l} (u,w)} < \infty $, where 
\begin{align*}
\|f\|_{M_{\alpha, \lambda}^{p, l}(u, w)}:=\sup 
_{r>0}\left[\int_{\mathbb{R}^{n}}\left(\varphi  
(Q(z,r))^{\alpha}u(Q(z,r))^{-\lambda}\left\|f\chi_{Q(z,r)}\right\|_{L^{p}(w)}\right)^{l} 
d z\right]^{1 / l}.
\end{align*}
The weak Morrey space $ WM_{\alpha ,\lambda }^{p,l} (u,w) $ is defined as the 
set of all measurable functions $ f $ satisfying $ \left \| f \right \| 
_{WM_{\alpha ,\lambda }^{p,l} (u,w)} < \infty $, where 
\begin{align*}
\|f\|_{WM_{\alpha, \lambda}^{p, l}(u, w)}:=\sup 
_{r>0}\left[\int_{\mathbb{R}^{n}}\left(\varphi  
(Q(z,r))^{\alpha}u(Q(z,r))^{-\lambda}\left\|f \chi_{Q(z,r)}\right\|_{L^{p,\infty 
}(w)}\right)^{l}dz\right]^{1 / l}.
\end{align*}
\end{definition}
When $ u=w $, we will denote $ M_{\alpha, \lambda}^{ p, l}(u, w) $ and $ WM_{\alpha, \lambda}^{ p, l}(u, w) $ by $M_{\alpha, \lambda}^{ p, l}(w) $ and $WM_{\alpha, \lambda}^{ p, l}(w)$ for brevity, respectively.

In 2021, Zhao and Zhou \cite{b38} studied new weighted norm inequalities for 
certain classes of multilinear operators on Morrey spaces. In this 
section, applying $ A_{\vec p}^{\infty } (\varphi )$, the boundedness on the 
new Morrey space is obtained for a new class of multilinear square 
operators with generalized kernels. The definition of the new Morrey space is 
given below.
\begin{definition}{\rm(\cite{b38})}
Let $ u $ be a weight, $ \vec w=(w_1,\ldots ,w_m) $ be a vector weight and $ 
\lambda \in [0,1) $, $ 1< l\le \infty  $, $ \alpha \in (- \infty ,\infty ) 
$, $ \vec p=(p_1, \ldots ,p_m) $, $ Q:=Q(z,r)$. The new strong Morrey space $ 
M_{\alpha ,\lambda }^{\vec p,l} (u,\vec w) $ is defined as the set of all  
vector-valued measurable functions $ \vec f=(f_1,\ldots ,f_m) $ satisfying $\| 
\vec f \| _{M_{\alpha ,\lambda }^{\vec p,l} (u,\vec w)} < \infty $, where 
\begin{align*}
\|\vec f\|_{M_{\alpha, \lambda}^{\vec p, l}(u, \vec w)}:=\sup 
_{r>0}\left[\int_{\mathbb{R}^{n}}\left(\varphi  
(Q(z,r))^{\alpha}u(Q(z,r))^{-\lambda} \prod_{j=1}^{m} \left\|f_j 
\chi_{Q(z,r)}\right\|_{L^{p_j}(w_j)}\right)^{l} d z\right]^{1 / l}.
\nonumber
\end{align*}

The new weak Morrey space $ WM_{\alpha ,\lambda }^{\vec p,l} (u,\vec w) $ is 
defined as the set of all vector-valued measurable functions $ \vec f=(f_1,\ldots 
,f_m) $ satisfying $\| \vec f  \| _{WM_{\alpha ,\lambda }^{\vec p,l} (u,\vec 
w)} < \infty$, where 
\begin{align*}
\|\vec f\|_{WM_{\alpha, \lambda}^{\vec p, l}(u, \vec w)}
&:=\sup _{r>0}\left[\int_{\mathbb{R}^{n}}\left(\varphi  
(Q(z,r))^{\alpha}u(Q(z,r))^{-\lambda} \prod_{j=1}^{m} \left\|f_j 
\chi_{Q(z,r)}\right\|_{L^{p_j,\infty }(w_j)}\right)^{l} d z\right]^{1/l }
\end{align*}
for $ p_j\ge 1 $, $ j=1,\ldots ,m $. When $ m=1 $, they are just the spaces $ 
M_{\alpha, \lambda}^{ p, l}(u, w) $ and $ WM_{\alpha, \lambda}^{ p, l}(u, w) $, 
respectively.
\end{definition}

\subsection{Necessary lemmas}

\begin{lemma}\label{Lemma3.1} {\rm(\cite{b38})}  
Let $ 0< \theta < \infty ,1\le p< \infty$ and $ E $ be any measurable subset  
of a cube $ Q $. If $ w\in A_{p}^{\theta } (\varphi ) $, then there exist 
positive constants $ 0< \delta < 1 $, $ \eta  $ and $ C $ such that 
\begin{align*}
\frac{w(E)}{w(Q)} \leq C 
\varphi(Q)^{\eta}\left(\frac{|E|}{|Q|}\right)^{\delta}.
\end{align*}
\end {lemma}
\begin{lemma}\label{Lemma3.2} {\rm(\cite{b38})}  
Let $ 0< \theta < \infty$, and $1\le p< \infty  $. If $ w\in A_{p}^{\theta 
}(\varphi ) $, then there exists a positive constant $ C $, such that for any 
$\rho > 1 $,
\begin{align*}
w(\rho Q) \leq C \varphi(\rho Q)^{p \theta} w(Q).
\end{align*}
\end {lemma}

\subsection{Main theorem and its proof}

\quad\quad  The boundedness on the new Morrey space is established for the 
new multilinear square operator with generalized kernel. We state our main 
results as follows.

\begin{theorem}\label{Theorem3.1} 
Suppose that $m\geq2$, $T$ is the multilinear square operator with generalized 
kernel as in Definition \ref{definition1.1} and  $ {\textstyle 
\sum_{k=1}^{\infty }}  C_k< \infty$. Let $ \vec p=(p_1,\ldots ,p_m) $, $  
\alpha \in (-\infty ,\infty )$, $ \lambda \in [0,1) $, $ 1\le p<l\le \infty $,  
$ 1/p = {\textstyle \sum_{j=1}^{m}}  1/p_j$, $ \vec w=(w_1,\ldots ,w_m)\in 
A_{\vec p /q'}^{\theta }(\varphi ) $, $ v_{\vec w}= {\textstyle 
\prod_{j=1}^{m}}  w_{j}^{p/p_j}$, and $\beta =\alpha +mp\theta \lambda/q'  $ 
with $ \theta > 0 $. Then the following statements hold.
\begin{itemize}
\item [\rm (i)] If $ q'< p_{i}< \infty  $, $ i=1,\ldots ,m $, then $ T $ is 
bounded 
from $ M_{\beta , \lambda}^{\vec p, l}\left ( v_{\vec{w}},\vec w \right )  $to 
$ M_{\alpha, \lambda}^{p, l}\left(v_{\vec{w}}\right) $.
	
\item [\rm (ii)] If $  q'\le  p_{i}< \infty  $, $ i=1,\ldots ,m $ and at least 
one of the $ p_i=q' $, then $ T $ is bounded from $ M_{\beta , \lambda}^{\vec 
p, l}\left ( v_{\vec{w}},\vec w \right )$ to $ WM_{\alpha, \lambda}^{p, 
l}\left(v_{\vec{w}}\right)$.
\end{itemize}
%
	\end {theorem}
	\setlength{\abovedisplayskip}{3pt}
	\setlength{\belowdisplayskip}{3pt}
	\noindent{\it \rm Proof.} (i) For any $ Q:=Q(z,r)\subset \mathbb{R}^n $, let $ 
	Q^\ast =3Q $. We split each $f_j$ as $f_j=f_j^0+f_j^\infty$, for 
	each $f_j^0=f_j \chi_{Q^*}$, $ f_j^\infty=f_j-f_j^0 $.
	\begin{align*}
		\prod_{i=1}^{m} f_{i}\left(y_{i}\right) 
		=\prod_{i=1}^m f_i^{0}(y_i)+\sum_{(\alpha_1,\ldots,\alpha_m)\in 
			\mathscr{L}}f_1^{\alpha_1}(y_1)\cdots f_m^{\alpha_m}(y_m),
	\end{align*}
	where  $\mathscr{L}:=\{(\alpha_1,\ldots,\alpha_m)$: there is at least one 
	$\alpha_i=\infty\}$. 
	
	We have for $x\in Q(z,r)$,
	\begin{align*}
		\left|T\left(f_{1},\ldots, f_{m}\right)(x)\right|
		&\leq \left|T\left(f_{1}^{0}, \ldots, f_{m}^{0}\right)(x)\right| +
		\sum_{(\alpha_1,\ldots,\alpha_m)\in \mathscr{L}} \left | 
		T\left(f_{1}^{\alpha_{1}}, 
		\ldots,f_{m}^{\alpha_{m}}\right)(x) \right | \\
		&:= \mathrm{I}+\mathrm{II}\\
		&:= \mathrm{I}+\sum_{(\alpha_1,\ldots,\alpha_m)\in \mathscr{L}}\mathrm{II}_{\alpha _1,\ldots,\alpha _m}.
	\end{align*}
	Taking $ L^{p} (\upsilon _{\vec \omega }) $ norm on the cube $ Q(z,r) $ of $ 
	\mathrm{I} $, by Theorem \ref{Theorem2.1} {\rm (i)}, we conclude 
	that
	\begin{align*}
		\left\|T\left(f_{1}^{0},\ldots, 
		f_{m}^{0}\right)\chi 
		_{Q(z,r)}\right\|_{L^{p}\left(v_{\vec{w}}\right)}\le 
		C\prod_{i=1}^{m}\left\| f_i\chi_{Q(z,3r)}\right\|_{L^{p_{i}}\left(w_{i}\right)}.
	\end{align*}
	
	Let us estimate the term $ \mathrm{II} $. Because each term of 
	$\mathscr{L}$ contains at least one $\alpha_i=\infty, i=1,\ldots,m$, we choose 
	one of them for discussion. Using Minkowski's inequality, the size 
	condition \eqref{1.2}, Hölder’s inequality and $ \vec w\in A_{\vec 
		p/q'}^{\theta  }(\varphi )$, we deduce that
	\begin{align}\label{3.1} 
		\mathrm{II} _{\alpha _1,\ldots,\alpha _m}
		&\le \left(\int_{0}^{\infty} \left ( \int_{(\mathbb{R}^{n})^{m}\backslash 
			(Q^\ast )^{m}} \left|K_{t}(x, \vec{y})\right|\prod_{j=1}^{m} 
		\left|f_{j}\left(y_{j}\right)\right| d \vec{y} \right ) ^{2} 
		\frac{d t}{t}\right)^{\frac{1}{2}} \nonumber\\
		&\leq \int_{(\mathbb{R}^{n})^{m}\backslash (Q^\ast)^{m}}\left(\int_{0}^{\infty} 
		\left |  K_{t}(x, \vec{y})\right|^{2} \frac{d 	t}{t}\right)^{\frac{1}{2} 
		} \prod_{j=1}^{m} \left|f_{j}\left(y_{j}\right)\right| d \vec{y} \nonumber 
		\\
	& \leq C \sum_{k=1}^{\infty } \int_{\left(3^{k+1}Q\right)^{m} 
		\backslash\left(3^k Q\right)^{m}}\frac{\left|f_{1}\left(y_{1}\right) \cdots 
		f_{m}\left(y_{m}\right)\right|}{\left(\sum_{i=1}^{m}\left|x-y_{i}\right|\right)^{m
			n}\left(1+\sum_{i=1}^{m}\left|x-y_{i}\right|\right)^{N}} d \vec{y} \nonumber\\
	& \leq C \sum_{k=1}^{\infty} \frac{1}{\left | 3^{k+1}Q \right |^{m}  
		\varphi\left(3^{k+1} Q\right)^{N}} \int_{\left(3^{k+1} Q\right)^{m}} 
	\prod_{i=1}^{m}\left|f_{i}\left(y_{i}\right)\right| d \vec{y} \\
	& \leq C\sum_{k=1}^{\infty} \varphi\left(3^{k+1} Q\right)^{-N}\left|3^{k+1} 
	Q\right|^{-\frac{m}{q'} }\prod_{i=1}^{m} \left ( \int_{3^{k+1} 
		Q}\left|f_{i}\left(y_{i}\right)\right|^{q'} d y_{i} \right )^{\frac{1}{q'} } 
	\nonumber \\
	& \leq C\sum_{k=1}^{\infty} \varphi\left(3^{k+1} Q\right)^{-N}\left|3^{k+1} 
	Q\right|^{-\frac{m}{q'} }
	\prod_{i=1}^{m} \left ( \int_{3^{k+1} 
		Q}\left|f_{i}\left(y_{i}\right)\right|^{p_i}w_i(y_i) d y_{i} \right 
	)^{\frac{1}{p_i} } \nonumber \\
	& \ \ \ \ \times\left ( \int_{3^{k+1} Q}w_i(y_i)^{1-\frac{p_i}{p_i-q'} } d 
	y_{i} \right )^{\frac{p_i-q'}{p_i} \cdot \frac{1}{q'} } \nonumber \\
	& \leq C\sum_{k=1}^{\infty} \frac{\varphi\left(3^{k+1} Q\right)^{\frac{m\theta 
			}{q'} -N}}{{\left(\int_{3^{k+1} Q} 
			v_{\vec{w}}\right)^{1 / p}}} \prod_{i=1}^{m}\left\|f_{i} \chi_{3^{k+1} 
		Q}\right\|_{L^{p_i} (w_i)}. \nonumber
\end{align}
Combining the estimates of the terms $ \mathrm{I} $ and $ \mathrm{II} $, we 
obtain 
\begin{equation}\label{3.2} 
	|T(\vec f)(x)|\le  | T(f_{1}^{0}, \ldots , f_{m}^{0})(x)| +C 
	\sum_{k=1}^{\infty} \frac{\varphi\left(3^{k+1} Q\right)^{{\frac{m 
					\theta}{q'} }- N}}{\left(\int_{3^{k+1} Q}	v_{\vec{w}}\right)^{1/ p}} 
	\prod_{i=1}^{m}\left\|f_{i} \chi_{3^{k+1} Q}\right\|_{L^{p_i} (w_i)}.
\end{equation}
Taking $ L^{p} (\upsilon _{\vec \omega }) $ norm on the cube $ Q(z,r) $ of 
\eqref{3.2}, we get
\begin{align}\label{3.3} 
	\left\|T(\vec{f}) \chi_{Q(z, r)}\right\|_{L^{p}\left(v_{\vec{w}}\right)} 
	&\leq C \prod_{i=1}^{m}\left\| 
	f_i\chi_{Q(z,3r)}\right\|_{L^{p_{i}}\left(w_{i}\right)} \nonumber \\ 
	&\quad+C \sum_{k=1}^{\infty} \frac{\varphi\left(3^{k+1} Q\right)^{{\frac{m 
					\theta}{q'} }- N}\cdot\left(\int_{Q}v_{\vec{w}}\right)^{1 / 
			p}}{\left(\int_{3^{k+1} Q}	v_{\vec{w}}\right)^{1/ p}} 
	\prod_{i=1}^{m}\left\|f_{i} \chi_{3^{k+1} Q}\right\|_{L^{p_i} (w_i)}.
\end{align}
Using Lemma \ref{Lemma2.3} and Lemma \ref{Lemma3.1}, 
\begin{align*}
	\frac{(\int_{Q} v_{\vec{w}})^{\frac{1}{p} }}{{\left(\int_{3^{k+1} Q} 
			v_{\vec{w}}\right)^{\frac{1}{p}}}} 
	=\frac{\big(v_{\vec{w}}(Q)\big) ^{\frac{1}{p} 
	}}{{\big(v_{\vec{w}}(3^{k+1}Q)\big)^{\frac{1}{p}}}} 
	\le C\frac{1}{3^{nk\delta /p}} \varphi (3^{k+1} Q)^{\frac{\eta }{p} }.
\end{align*}
Using Lemma \ref{Lemma2.3} and Lemma \ref{Lemma3.2}, $ v_{\vec w}\in 
A_{mp/q'}^{\theta } (\varphi ) $, and we conclude
\begin{align*}
	v_{\vec w}(Q)^{-\lambda }\le C\varphi (3^{k+1}Q)^{\frac{mp\theta \lambda }{q'}  
	}v_{\vec w}(3^{k+1}Q)^{-\lambda }.
\end{align*}
Noting that $ \alpha \in (-\infty  ,\infty )  $, it follows that
$$\varphi (Q)^{\alpha }\le \varphi (3^{k+1}Q)^{|\alpha| },$$
and
\begin{align*}
	\frac{\varphi (Q)^{\alpha }}{ \varphi (3^{k+1}Q)^{\frac{N}{2}  }} \le 
	\frac{\varphi (3^{k+1}Q)^{|\alpha| }}{\varphi (3^{k+1}Q)^{\frac{N}{2}  }}
	=\frac{\varphi (3^{k+1}Q)^{\alpha }}{\varphi (3^{k+1}Q)^{\frac{N}{2} 
			-|\alpha|+\alpha  }}.
\end{align*}
Taking $ N\ge 2(|\alpha |-\alpha ) $, we obtain
$$\frac{\varphi (Q)^{\alpha }}{ \varphi (3^{k+1}Q)^{\frac{N}{2}  }} \le C 
\varphi (3^{k+1}Q) ^{\alpha}.$$
Multiplying both sides of \eqref{3.3} by $ \varphi(Q)^{\alpha} 
v_{\vec{w}}(Q)^{-\lambda} $, by Lemma \ref{Lemma2.3}, Lemma \ref{Lemma3.1} and 
Lemma \ref{Lemma3.2}, and taking $ N\ge \mathop{\max}\left \{ 2(|\alpha 
|-\alpha ), 2m\theta/q'+2\eta/p \right \}$, we obtain
\begin{align}\label{3.4} 
	&\varphi(Q)^{\alpha} v_{\vec{w}}(Q)^{-\lambda}\left\|T(\vec{f}) \chi_{Q(z, 
		r)}\right\|_{L^{p}\left(v_{\vec{w}}\right)} \nonumber\\
	& \leq C\varphi (Q)^{\alpha }v_{\vec{w}}(Q)^{-\lambda}  
	\prod_{i=1}^{m}\left\|f_{i} 
	\chi_{3Q}\right\|_{L^{p_{i}}\left(w_{i}\right)} \nonumber\\
	& \quad +C \sum_{k=1}^{\infty} \frac{\varphi\left(3^{k+1} Q\right)^{{\frac{m 
					\theta}{q'} }- N}\cdot\left(\int_{Q}v_{\vec{w}}\right)^{1 / 
			p}}{\left(\int_{3^{k+1} Q}v_{\vec{w}}\right)^{1/ p}} 
	\varphi(Q)^{\alpha} 
	v_{\vec{w}}(Q)^{-\lambda}
	\prod_{i=1}^{m}\left\|f_{i} \chi_{3^{k+1} Q}\right\|_{L^{p_i} (w_i)}\nonumber\\
	&\leq C\varphi (Q)^{\alpha }\varphi\left(3 Q\right)^{\frac{mp\theta \lambda 
		}{q'} } v_{\vec{w}}\left(3 Q\right)^{-\lambda} \prod_{i=1}^{m}\left\|f_{i} 
	\chi_{3Q}\right\|_{L^{p_{i}}\left(w_{i}\right)} \nonumber\\
	& \quad + C \sum_{k=1}^{\infty} \frac{1}{3^{n k \delta / p}} 
	\varphi\left(3^{k+1} Q\right)^{\frac{m\theta }{q'}-\frac{N}{2}+\frac{\eta 
		}{p}   }
	\frac{\varphi (Q)^{\alpha }}{\varphi\left(3^{k+1} Q\right)^{\frac{N}{2} }} 
	\varphi\left(3^{k+1} Q\right)^{\frac{mp\theta \lambda }{q'} } 
	v_{\vec{w}}\left(3^{k+1} Q\right)^{-\lambda} \nonumber\\
	&\quad \times \prod_{i=1}^{m}\left\|f_{i} 
	\chi_{3^{k+1}Q}\right\|_{L^{p_{i}}\left(w_{i}\right)} \nonumber\\
	&\leq C\varphi\left(3 Q\right)^{\alpha+\frac{mp\theta \lambda }{q'} } 
	v_{\vec{w}}\left(3 Q\right)^{-\lambda} \prod_{i=1}^{m}\left\|f_{i} 
	\chi_{3Q}\right\|_{L^{p_{i}}\left(w_{i}\right)} \nonumber\\
	&\quad + C \sum_{k=1}^{\infty} \frac{1}{3^{n k \delta / p}} 
	\varphi\left(3^{k+1} Q\right)^{\alpha+\frac{mp\theta \lambda }{q'} } 
	v_{\vec{w}}\left(3^{k+1} Q\right)^{-\lambda} \prod_{i=1}^{m}\left\|f_{i} 
	\chi_{3^{k+1}Q}\right\|_{L^{p_{i}}\left(w_{i}\right)} \nonumber\\
	& \leq C \sum_{k=0}^{\infty} \frac{1}{3^{n k \delta / p}} 
	\varphi\left(3^{k+1} Q\right)^{\alpha+\frac{mp\theta \lambda }{q'} } 
	v_{\vec{w}}\left(3^{k+1} Q\right)^{-\lambda} \prod_{i=1}^{m}\left\|f_{i} 
	\chi_{3^{k+1}Q}\right\|_{L^{p_{i}}\left(w_{i}\right)}.
\end{align}
Taking $L^{l}(\mathbb{R}^{n})$ norm of both sides of \eqref{3.4}, we have
\begin{align*}
	&\left\|\varphi(Q)^{\alpha} v_{\vec{w}}(Q)^{-\lambda}|| T(\vec{f}) 
	\chi_{Q(z,r)}||_{L^{p}\left(v_{\vec{w}}\right)}\right\|_{L^{l}\left
		(\mathbb{R}^{n}\right)} \\
	& \leq C \sum_{k=0}^{\infty} \frac{1}{3^{n k \delta/ p}} 
	\left\|\varphi\left(3^{k+1} Q\right)^{\alpha+\frac{mp\theta 
			\lambda }{q'} } v_{\vec{w}}\left(3^{k+1} Q\right)^{-\lambda  }\prod_{i=1}^{m}|| f_{i} \chi_{3^{k+1} 
		Q}||_{L^{p_{i}}\left(w_{i}\right)}\right\|_{L^{l}\left(\mathbb{R}^{n}\right)}.
\end{align*}
Since $ \beta =\alpha +mp\theta \lambda /q'$, then Theorem 
\ref{Theorem3.1} {\rm (i)} is obtained.

\noindent{\it \rm (ii)} For any $ t > 0$, it 
follows that
\begin{align*}
	v_{\vec{w}}(&\{x \in Q(z, r):|T(\vec{f})(x)|>t\})^{1 / p} \\
	& \leq v_{\vec{w}}\left(\left\{x \in Q(z, r):\left|T\left(f_{1}^{0}, \ldots , 
	f_{m}^{0}\right)(x)\right|>t / 2^{m}\right\}\right)^{1 / p} \\
	&\ \ \ \ +\sum_{(\alpha_1,\ldots,\alpha_m)\in \mathscr{L}} 
	v_{\vec{w}}\left(\left\{x \in Q(z, r):\left|T\left(f_{1}^{\alpha_{1}},\ldots , 
	f_{m}^{\alpha_{m}}\right)(x)\right|>t / 2^{m}\right\}\right)^{1 / p}\\
	&:= \mathrm{III}+\mathrm{IV} ,
\end{align*}
where  $\mathscr{L}:=\{(\alpha_1,\ldots,\alpha_m)$: there is at least one 
$\alpha_i=\infty\}$. 

For the term $ \mathrm{III} $, by Theorem \ref{Theorem2.1} {\rm (ii)}, we obtain
\begin{align*}
	\mathrm{III} \leq \frac{C}{t} 
	\prod_{i=1}^{m}\left(\int_{Q^{*}}\left|f_{i}(x)\right|^{p_{i}} w_{i}(x) d 
	x\right)^{1 / p_{i}}=\frac{C}{t} \prod_{i=1}^{m}\left\|f_{i} \chi_{ 
		3Q}\right\|_{L^{p_{i}}\left(w_{i}\right)}.
\end{align*}

For the term $ \mathrm{IV} $, the following pointwise estimate can be 
obtained from \eqref{3.1}.
\begin{align*}
	\left |T\left(f_{1}^{\alpha_{1}}, \ldots ,f_{m}^{\alpha_{m}}\right)(x) \right 
	| \leq C \sum_{k=1}^{\infty} \frac{1}{  
		\varphi\left(3^{k+1} Q\right)^{N}} \prod_{i=1}^{m}\frac{1}{\left | 3^{k+1}Q 
		\right |} \int_{3^{k+1} Q} 
	\left|f_{i}\left(y_{i}\right)\right| d {y_i}.
\end{align*}

If $ p_i> q' $, it follows that
{\setlength\abovedisplayskip{10pt}
	\setlength\belowdisplayskip{10pt}
	\begin{align*}
		&\frac{1}{\left | 3^{k+1}Q \right |} \int_{3^{k+1} Q} 
		\left|f_{i}\left(y_{i}\right)\right| d {y_i} \\
		&\le \left ( \frac{1}{\left | 3^{k+1}Q \right |} \int_{3^{k+1} Q} 
		\left|f_{i}\left(y_{i}\right)\right|^{q'} d {y_i} \right ) ^{\frac{1}{q'} }\\
		&\le |3^{k+1}Q|^{-\frac{1}{q'} 
		}\left(\int_{3^{k+1}Q}\left|f_{i}\left(y_{i}\right)\right| ^{p_{i}}
		w_{i}\left(y_{i}\right) d y_{i}\right)^{\frac{1}{p_{i}} 
		}\left(\int_{3^{k+1} Q} 
		w_{i}\left(y_{i}\right)^{1-(\frac{p_i}{q'} 
			)^{'}}dy_{i}\right)^{\frac{1}{(\frac{p_i}{q'} )^{'}}\cdot \frac{1}{q'}  }\\
		&=\left(\frac{1}{\varphi\left(3^{k+1} 
			Q\right)^{\theta }|3^{k+1}Q|} \int_{3^{k+1} Q} 
		w_{i}\left(y_{i}\right)^{1-(\frac{p_i}{q'} 
			)^{'}}dy_{i}\right)^{\frac{1}{(\frac{p_i}{q'} )^{'}}\cdot \frac{1}{q'}  }\\
		&\ \ \ \ \times \varphi\left(3^{k+1} 
		Q\right)^{\theta(\frac{1}{q'}-\frac{1}{p_i}  ) }|3^{k+1}Q|^{-\frac{1}{p_i} }
		\left\|f_{i} \chi_{3^{k+1} 
			Q}\right\|_{L^{p_{i}}\left(w_{i}\right)}.
\end{align*}

If $ p_i= q' $, then
\begin{align*}
&\frac{1}{\left | 3^{k+1}Q \right |} \int_{3^{k+1} Q} 
\left|f_{i}\left(y_{i}\right)\right| d {y_i} \\
&=\frac{1}{\left | 3^{k+1}Q \right |} \int_{3^{k+1} Q} 
\left|f_{i}\left(y_{i}\right)\right|w_i(y_i)^{\frac{1}{p_i} } 
w_i(y_i)^{-\frac{1}{p_i} }d {y_i}\\
&\le \left ( \frac{1}{\left | 3^{k+1}Q \right |} \int_{3^{k+1} Q} 
\left|f_{i}\left(y_{i}\right)\right|w_i(y_i)^{\frac{1}{p_i} } d {y_i} \right ) 
 \left ( \inf_{y_i\in3^{k+1} Q }w_i(y_i)\right )^{-\frac{1}{p_i} 
} \\
&\le \left ( \frac{1}{\left | 3^{k+1}Q \right |} \int_{3^{k+1} Q} 
\left|f_{i}\left(y_{i}\right)\right|^{p_i}w_i(y_i) d {y_i} \right 
)^{\frac{1}{p_i} }  \left ( \inf_{y_i\in3^{k+1} Q }w_i(y_i)\right 
)^{-\frac{1}{p_i} } \\
&=\left\|f_{i} \chi_{3^{k+1} 
Q}\right\|_{L^{p_{i}}\left(w_{i}\right)} \left [ \left ( \inf_{y_i\in3^{k+1} Q 
}w_i(y_i)\right )^{-1 }  \right ]^{\frac{1}{q'} }|3^{k+1} Q| ^{-\frac{1}{p_i} }.
\end{align*}
To sum up, we have
\begin{align*}
\left |T\left(f_{1}^{\alpha_{1}}, \ldots ,f_{m}^{\alpha_{m}}\right)(x) \right |
&\le C \sum_{k=1}^{\infty} \frac{1}{  
\varphi\left(3^{k+1} Q\right)^{N}} \\
&\ \ \ \ \times \left ( \prod_{i=1}^{m}\left\|f_{i} 
\chi_{3^{k+1} Q}\right\|_{L^{p_{i}}\left(w_{i}\right)}\varphi\left(3^{k+1} 
Q\right)^{\theta(\frac{1}{q'}-\frac{1}{p_i}  ) }|3^{k+1}Q|^{-\frac{1}{p_i} }  
\right ) \\
&\ \ \ \ \times \prod_{i=1}^{m}  \left(\frac{1}{\varphi\left(3^{k+1} 
Q\right)^{\theta }|3^{k+1}Q|} \int_{3^{k+1} Q} 
w_{i}\left(y_{i}\right)^{1-(\frac{p_i}{q'} 
)^{'}}dy_{i}\right)^{\frac{1}{(\frac{p_i}{q'} )^{'}}\cdot \frac{1}{q'}  },
\end{align*}
where the term $$\left ( \frac{1}{|3^{k+1}Q|} \int_{3^{k+1} Q} 
w_{i}\left(y_{i}\right)^{1-(\frac{p_i}{q'} 
)^{'}}dy_{i} \right ) ^{\frac{1}{(\frac{p_i}{q'} )^{'}} }$$ is considered to be 
$\left ( \inf_{y_i\in3^{k+1} Q }w_i(y_i)\right )^{-1 } $ when $ p_i=q' $.

Thus, by the definition of $  \vec w\in A_{\vec p/q'}^{\theta}(\varphi ) $, 
there is
\begin{align*}
\left |T\left(f_{1}^{\alpha_{1}},\ldots,f_{m}^{\alpha_{m}}\right)(x) \right |
& \le C \sum_{k=1}^{\infty} \frac{1}{\varphi\left(3^{k+1} Q\right)^{N}} \\
&\ \ \ \ \times \left ( \prod_{i=1}^{m}\left\|f_{i} 
\chi_{3^{k+1} Q}\right\|_{L^{p_{i}}\left(w_{i}\right)}\varphi\left(3^{k+1} 
Q\right)^{\theta(\frac{1}{q'}-\frac{1}{p_i}  ) }|3^{k+1}Q|^{-\frac{1}{p_i} }  
\right ) \\
&\ \ \ \ \times  \left({\varphi\left(3^{k+1} 
Q\right)^{\theta }|3^{k+1}Q|} v_{\vec w}(3^{k+1} 
Q)^{-1}\right)^{\frac{q'}{p}\cdot \frac{1}{q'}  }\\
&=C \sum_{k=1}^{\infty} \frac{\varphi\left(3^{k+1} Q\right)^{{\frac{m 
		\theta}{q'} }- N}}{\left(\int_{3^{k+1} Q} v_{\vec{w}}\right)^{1 / p}} 
\prod_{i=1}^{m}\left\|f_{i} \chi_{3^{k+1} 
Q}\right\|_{L^{p_{i}}\left(w_{i}\right)}.
\end{align*}
Therefore, we get
\begin{align*}
v_{\vec{w}} &\left(\left\{x \in Q(z, r):\left|T\left(f_{1}^{\alpha_{1}}, 
\ldots , f_{m}^{\alpha_{m}}\right)(x)\right|>t / 2^{m}\right\}\right)^{1 
/p} \\
&=\left(\int_{\left\{x \in Q(z, r):\left|T\left(f_{1}^{\alpha_{1}},\ldots , 
f_{m}^{\alpha_{m}}\right)(x)\right|>t / 2^{m}\right\}} v_{\vec{w}}(x) d 
x\right)^{1 / p} \\
& \leq C\left(\int_{Q(z, r)} \frac{\left|T\left(f_{1}^{\alpha_{1}},\ldots , 
f_{m}^{\alpha_{m}}\right)(x)\right|^{p}}{t^{p}} v_{\vec{w}}(x) d 
x\right)^{1 / p} \\
& \leq C\left(t^{-p} \int_{Q(z, r)}\left(\sum_{k=1}^{\infty} 
\frac{\varphi\left(3^{k+1} Q\right)^{{\frac{m \theta}{q'} }- N}
}{\left(\int_{3^{k+1} Q} v_{\vec{w}}\right)^{1 / p}} 
\prod_{i=1}^{m}\left\|f_{i} \chi_{3^{k+1} 
Q}\right\|_{L^{p_{i}}(w_{i})}\right)^{p} v_{\vec{w}}(x) d 
x\right)^{1 / p} \\
& \leq C t^{-1} \sum_{k=1}^{\infty} \frac{\varphi\left(3^{k+1} 
Q\right)^{{\frac{m \theta}{q'}}- N}\cdot \left(\int_{Q} 
v_{\vec{w}}\right)^{1 / p}}{\left(\int_{3^{k+1} Q} v_{\vec{w}}\right)^{1 / p}} 
\prod_{i=1}^{m}\left\|f_{i} \chi_{3^{k+1} 
Q}\right\|_{L^{p_{i}}\left(w_{i}\right)}.
\end{align*}
Summing up the above estimates for the terms $\mathrm{III}$ and $\mathrm{IV} 
$, it yields that
\begin{align}\label{3.5} 
&\left\|T(\vec{f}) \chi_{Q(z, r)}\right\|_{L^{p, 
\infty}\left(v_{\vec{w}}\right)} \nonumber\\
&=\sup_{t > 0} t v_{\vec{w}}(\{x \in Q(z, r):|T(\vec{f})(x)|>t\})^{1 / p} 
\nonumber\\
&\leq C \prod_{i=1}^{m}\left\| 
f_i\chi_{3Q}\right\|_{L^{p_{i}}\left(w_{i}\right)} \nonumber\\
&\ \ \ \ + C \sum_{k=1}^{\infty} \frac{\varphi\left(3^{k+1} Q\right)^{{\frac{m 
		\theta}{q'}}- N}\cdot\left(\int_{Q} v_{\vec{w}}\right)^{1 / 
p}}{\left(\int_{3^{k+1} Q} v_{\vec{w}}\right)^{1 / p}} 
\prod_{i=1}^{m}\left\|f_{i}\chi_{3^{k+1}Q}\right\|_{L^{p_{i}}\left(w_{i}\right)}.
\end{align}
Multiplying both sides of \eqref{3.5} by $ \varphi(Q)^{\alpha} 
v_{\vec{w}}(Q)^{-\lambda} $, the following estimate can be obtained
similarly to \eqref{3.4}.
\begin{align}\label{3.6}
	&\varphi(Q)^{\alpha} v_{\vec{w}}(Q)^{-\lambda}\left\|T(\vec{f}) \chi_{Q(z, 
		r)}\right\|_{L^{p,\infty}\left(v_{\vec{w}}\right)} \nonumber\\
	& \leq C \sum_{k=0}^{\infty} \frac{1}{3^{n k \delta / p}} 
	\varphi\left(3^{k+1} Q\right)^{\alpha+\frac{mp\theta \lambda }{q'} } 
	v_{\vec{w}}\left(3^{k+1} Q\right)^{-\lambda} \prod_{i=1}^{m}\left\|f_{i} 
	\chi_{3^{k+1}Q}\right\|_{L^{p_{i}}\left(w_{i}\right)}.
\end{align}
Taking $L^{l}(\mathbb{R}^{n})$ norm of both sides of \eqref{3.6}, we get
\begin{align*}
	&\left\|\varphi(Q)^{\alpha} v_{\vec{w}}(Q)^{-\lambda}|| T(\vec{f}) 
	\chi_{Q(z,r)}||_{L^{p,\infty}\left(v_{\vec{w}}\right)}\right\|_{L^{l}\left
		(\mathbb{R}^{n}\right)} \\
	& \leq C \sum_{k=0}^{\infty} \frac{1}{3^{n k \delta/ p}} 
	\left\|\varphi\left(3^{k+1} Q\right)^{\alpha+\frac{mp\theta 
			\lambda }{q'} } v_{\vec{w}}\left(3^{k+1} Q\right)^{-\lambda  }\prod_{i=1}^{m}|| f_{i} \chi_{3^{k+1} 
		Q}||_{L^{p_{i}}\left(w_{i}\right)}\right\|_{L^{l}\left(\mathbb{R}^{n}\right)}.
\end{align*}

The proof of Theorem \ref{Theorem3.1} is finished.

\bigskip

\noindent Chunliang Li,   Shuhui Yang   and  Yan Lin
\smallskip

\noindent School of Science, China University of Mining and
Technology, Beijing 100083,  People's Republic of China

\smallskip

\noindent{\it E-mails:} \texttt{lichunliang@student.cumtb.edu.cn} (C. Li)

\noindent\phantom{{\it E-mails:} }\texttt{yangshuhui@student.cumtb.edu.cn} (S. 
Yang)

\noindent\phantom{{\it E-mails:} }\texttt{linyan@cumtb.edu.cn} (Y. Lin)

\end{document}